\documentclass[12pt]{article}
\usepackage{caption}
\usepackage{amsmath}
\usepackage{amssymb}
\usepackage{ifthen}
\usepackage{bm}
\usepackage{amsfonts}
\usepackage{multirow}
\usepackage{latexsym}
\usepackage{amsthm}
\usepackage{mathrsfs}
\usepackage{color}
\usepackage[dvips]{graphicx}
\usepackage{float}
\usepackage{varioref}
\usepackage[top=1in, bottom=1in, left=1in, right=1in]{geometry}
\usepackage{array}
\newcommand{\Dchaintwo}[4]{
\rule[-3\unitlength]{0pt}{8\unitlength}
\begin{picture}(14,5)(0,3)
\put(1,2){\ifthenelse{\equal{#1}{l}}{\circle*{2}}{\circle{2}}}
\put(2,2){\line(1,0){10}}
\put(13,2){\ifthenelse{\equal{#1}{r}}{\circle*{2}}{\circle{2}}}
\put(1,5){\makebox[0pt]{\scriptsize #2}}
\put(7,4){\makebox[0pt]{\scriptsize #3}}
\put(13,5){\makebox[0pt]{\scriptsize #4}}
\end{picture}}

\newcommand{\Dchainthree}[6]{
\rule[-3\unitlength]{0pt}{8\unitlength}
\begin{picture}(26,5)(0,3)
\put(1,2){\ifthenelse{\equal{#1}{l}}{\circle*{2}}{\circle{2}}}
\put(2,2){\line(1,0){10}}
\put(13,2){\ifthenelse{\equal{#1}{m}}{\circle*{2}}{\circle{2}}}
\put(14,2){\line(1,0){10}}
\put(25,2){\ifthenelse{\equal{#1}{r}}{\circle*{2}}{\circle{2}}}
\put(1,5){\makebox[0pt]{\scriptsize #2}}
\put(7,4){\makebox[0pt]{\scriptsize #3}}
\put(13,5){\makebox[0pt]{\scriptsize #4}}
\put(19,4){\makebox[0pt]{\scriptsize #5}}
\put(25,5){\makebox[0pt]{\scriptsize #6}}
\end{picture}}
\newcommand{\Dtriangle}[7]{
\rule[-3\unitlength]{0pt}{12\unitlength}
\begin{picture}(18,7)(0,3)
\put(4,4){\ifthenelse{\equal{#1}{l}}{\circle*{2}}{\circle{2}}}
\put(5,4){\line(1,0){8}}
\put(14,4){\ifthenelse{\equal{#1}{r}}{\circle*{2}}{\circle{2}}}
\put(4.4472,4.8944){\line(1,2){4.1056}}
\put(9,14){\ifthenelse{\equal{#1}{t}}{\circle*{2}}{\circle{2}}}
\put(13.5528,4.8944){\line(-1,2){4.1056}}
\put(2,3){\makebox[0pt][r]{\scriptsize #2}}
\put(9,17){\makebox[0pt]{\scriptsize #3}}
\put(16,3){\makebox[0pt][l]{\scriptsize #4}}
\put(6,9){\makebox[0pt][r]{\scriptsize #5}}
\put(12.5,9){\makebox[0pt][l]{\scriptsize #6}}
\put(9,1){\makebox[0pt]{\scriptsize #7}}
\end{picture}}

\newcounter{cthm}%

\newlength{\mpb}

\newcommand{\Aut }{\mathrm{Aut}}
\newcommand{\Hom }{\mathrm{Hom}}
\newcommand{\End }{\mathrm{End}}

\newcommand{\cB }{\mathcal{B}}
\newcommand{\cC }{\mathcal{C}}
\newcommand{\cR }{\mathcal{R}}

\newcommand{\cD }{\mathcal{D}}

\newcommand{\cX }{\mathcal{X}}
\newcommand{\cY }{\mathcal{Y}}
\newcommand{\cW }{\mathcal{W}}

\newcommand{\lact }{.}

\newcommand{\lcoa }{\delta }

\newcommand{\ndN }{\mathbb{N}}

\newcommand{\ndZ }{\mathbb{Z}}

\newcommand{\ot }{\otimes }

\newcommand{\PBW }{Poincar\'e--Birkhoff--Witt }

\newcommand{\roots }{\boldsymbol{\Delta }}

\newcommand{\YD }{Yetter--Drinfel'd }

\newcommand{\ydD }{ {}^{G}_{G}\mathcal{YD}}

\newcommand{\id}{\mathrm{id}}
\newcommand{\al }{\alpha }

\newcommand{\btxandshort}[1]{and}%
\newcommand{\btxpagesshort}[1]{pp.}%
\newcommand{\Btxinshort}[1]{In}%
\newcommand{\btxphdthesis}[1]{phd-thesis}%
\newcommand{\btxeditorshort}[1]{Ed.}%
\newcommand{\btxeditorsshort}[1]{Eds.}%
\newcommand{\btxvolumeshort}[1]{vol.}%
\newcommand{\btxofseriesshort}[1]{ser.}%
\newcommand{\ffg }{\mathcal{F}_\theta^{G} }   
\newcommand{\fiso }{\mathcal{X}_\theta }      
\newcommand{\rersys }[1]{\boldsymbol{\Delta }{}^{#1\,\mathrm{re}}}
\newcommand{\rsys }{\boldsymbol{\Delta }}
\newcommand{\ad }{\mathrm{ad}}
\newcommand{\hSL}[1]{\widehat{s\ell}(#1)}
\newtheorem{theorem}{\bf Theorem}[section]
\newtheorem{lemma}[theorem]{\bf Lemma}
\newtheorem{prop}[theorem]{\bf Proposition}
\newtheorem{coro}[theorem]{\bf Corollary}

\newtheorem{defn}[theorem]{\bf Definition}

\newtheorem{remark}[theorem]{\bf Remark}

\begin{document}

\title{Rank three Nichols algebras of diagonal type over arbitrary fields}

\author{Jing Wang\thanks{supported by China Scholarship Council}
\\jing@mathematik.uni-marburg.de
\\Philipps-Universit\"at Marburg \\FB Mathematik und Informatik\\Hans-Meerwein-Stra\ss e\\35032 Marburg, Germany}
\date{}
\maketitle
\begin{abstract}
Over fields of arbitrary characteristic we classify all rank three Nichols algebras of diagonal type with a finite root system.
Our proof uses the classification of the finite Weyl groupoids of rank three.

\end{abstract}

\section*{Introduction}
The theory of Nichols algebras is motivated by the Hopf algebra theory and plays an important role in quantum groups~\cite{inp-Andr02,a-AndrGr99,inp-AndrSchn02, a-Schauen96}.
The structure of Nichols algebras was first introduced systematically by W.~Nichols~\cite{n-78} in 1978,
where he studied Hopf algebras.
Nichols algebras have been redescribed independently by S.L.~Woronowicz, M.~Rosso, S.~Majid and G.~Lusztig in many different ways,
see for example~\cite{w1987,w1989,Rosso98,maj05,l-2010}.
Besides that,
Nichols algebras have interesting applications to other research fields such as Kac-Moody Lie superalgebras~\cite[Example~3.2]{inp-Andr14} and conformal field theory~\cite{Semi-2011,Semi-2012,Semi-2013}.
Nichols algebras appeared naturally in the classification of pointed Hopf algebras by the lifting method of N.~Andruskiewitsch and H.-J.~Schneider~\cite{a-AndrSchn98,inp-AndrSchn02}.
The crucial step to classify pointed Hopf algebras is to determine all braided vector space $V$ such that the Nichols algebra $\cB(V)$ is finite dimensional.
Several authors obtained the classification result for infinite and finite dimensional Nichols algebra of Cartan type,
see~\cite{a-AndrSchn00,a-Heck06a,Rosso98}.
I.~Heckenberger classified all finite dimensional Nichols algebra of diagonal type in a series of papers~\cite{a-Heck04aa, a-Heck04bb,a-Heck04d, a-Heck09}.
The explicit defining generators and relations of such Nichols algebras were given~\cite{Ang1,Ang2}.
With the classification result~\cite{a-Heck09},
N.~Andruskiewitsch and H.-J.~Schneider~\cite{a-AndrSchn05} obtained a classification result about finite-dimensional pointed Hopf algebras under some technical assumptions.
Based on such successful applications,
the analyze to Nichols algebras over arbitrary fields is crucial and has also potential applications.
Towards this direction,
the authors~\cite{clw} discovered a combinatorial formula to study the relations in Nichols algebras and found new examples of Nichols algebras.
Over fields of arbitrary characteristic,
a complete list of rank two finite dimensional Nichols algebras of diagonal type were determined in~\cite{WH-14}.
In this paper,
we give the complete classification result of rank three case.

The important theoretical tools for the classification of Nichols algebra $\cB(V)$ are the root system and the Weyl groupoid of $\cB(V)$.
I.~Heckenberger~\cite{a-Heck06a} defined the root system and Weyl groupoid of $\cB(V)$ of diagonal type
which based on the ideas of V.~Kharchenko on \PBW basis of $\cB(V)$ of diagonal type~\cite[Theosrem~2]{a-Khar99}.
Later,
the combinatorial theory of
these two structures was initiated in~\cite{c-Heck09b,Y-Heck08a} by M.~Cuntz, I.~Heckenberger and H.~Yamane.
The theory of
root systems and Weyl groupoids was then carried out in more general Nichols
algebras~\cite{a-AHS08,HS10,a-HeckSchn12a}.
Further,
all finite Weyl groupoids were classified
in~\cite{c-Heck12a,c-Heck14a}.

We organize the paper as follows.
In section~\ref{se:Pre} the main notations and some general results are recalled.
One general result implies that we can attach a connected finite Cartan graph of rank $\theta$ to a tuple $M=(M_1, \dots, M_{\theta})$ of
one-dimensional \YD modules over an abelian group,
see Section~\ref{se:Pre} for the definitions.
If $\theta= 1$,
then the Cartan graph contains no information about $\cB(M)$.
Therefore we give the restriction $\theta\geq 2$.
Since the classification of $\theta=2$ was performed in~\cite{WH-14},
we record that the we add indecomposability assumption on the
semi-Cartan graph.
The indecomposability allows us to exclude
the components of the Cartan graph of lower rank.
In order to explicitly characterize finite connected indecomposable Cartan graphs of rank three,
we introduce three concepts,
good $A_3$ neighborhood,
good $B_3$ neighborhood,
and good $C_3$ neighborhood,
see Definitions~\ref{defA3},~\ref{defB3} and~\ref{defC3}.
In Proposition~\ref{prop:A3},
we prove that every finite connected indeconposable Cartan graph of rank three contains a point which has at least one of the good neighborhoods.
In Section~\ref{se:clasi} we formulate the main classification Theorem~\ref{Theo:clasi} and give a complete list of the Dynkin diagrams for rank three braided vector spaces of diagonal type with a finite root system over arbitrary fields,
see Tables~(\ref{tab.1})- (\ref{tab.3}).
Proposition~\ref{prop:A3} allows us to avoid complicated computations in the final proof of Theorem~\ref{Theo:clasi}.
The proof of our main classification theorem uses the classification result on the rank three Weyl groupoids~\cite{c-Heck12a}.
As a corollary of Theorem~\ref{Theo:clasi},
all finite
dimensional rank three Nichols algebras of diagonal type over a field of
positive characteristic are given,
see Corollary~\ref{coro-cla}.

Throughout this paper $\Bbbk$ always denotes a field of characteristic $p> 0$.
Let $\Bbbk^*=\Bbbk\setminus \{0\}$.
The set of natural numbers strictly greater than $0$ is denoted by $\ndN$ and then write $\ndN_0:=\ndN\bigcup \{0\}$.
For $n\in \ndN$,
the symbol $G_n'$ denotes the set of primitive $n$-th roots of unity
in $\Bbbk$,
that is $G'_n=\{q\in \Bbbk^*|\,\, q^n=1, q^k\not=1~\text{for all}~ 1\leq k < n\}.$
For $q\in \ndN$,
the symbol $(n)_q$ denotes the integer $1+q+\cdots+q^{n-1}$ for any $n\in \ndN$,
which is $0$ if and only if $q^n=1$ for $q\not=1$ or $p|n$ for $q=1$.
Notice that $(1)_q\not=0$ for any $q\in \ndN$.

The author would like to gratefully thank Professor I.~Heckenberger for the manuscript correction and providing valuable advices, as well as his support and encouragement during the author's pursuit of the PhD degree.

\section{Preliminaries}\label{se:Pre}
\subsection{Cartan graphs, root systems, and Weyl groupoids}\label{se:Cschemes}
We start by recalling the notations of semi-Cartan graphs,
root systems and Weyl groupoids,
mainly following the terminology introduced in~\cite{WH-14}.

Let $\theta\in \ndN$ and $I=\{1, \dots, \theta\}$.
We say that a generalized Cartan matrix $A\in \ndZ ^{I \times I}$ is \textbf{decomposable} if there exists a nonempty proper subset $I_1\subset I$ such that $a_{ij}=0$ whenever $i \in I_1$ and $j \notin I_1$.
Let $\cX$ be a non-empty set and $A^X =(a_{ij}^X)_{i,j \in I}$ a generalized Cartan
  matrix for any $X\in \cX$.
  For any $i\in I$,
  let $r_i : \cX \to \cX$,
  $ X \mapsto r(i,X)$,
  where $r: I\times \cX \to \cX$ is a map.
  The quadruple
  $$\cC = \cC (I, \cX, r, (A^X)_{X \in \cX})$$
  is called a \textbf{semi-Cartan graph} if
  $r_i^2 = \id_{\cX}$ for all $i \in I$,
  and $a^X_{ij} = a^{r_i(X)}_{ij}$ for all $X\in \cX$ and
  $i,j\in I$.
  We say that a semi-Cartan graph is~\textbf{indecomposable}
  if there exists $X \in \cX$ such that $A^X$ is not decomposable.

For the remaining part of this section,
let $\cC = \cC (I, \cX, r, (A^X)_{X \in \cX})$ be an indecomposable semi-Cartan graph.
For any point $X\in \cX$,
the elements of the set $\{r_i(X), i\in I\}$ are the \textbf{neighbors} of $X$.
The cardinality of $I$ is the \textbf{rank} of $\cC$ and the elements of $\cX$ are called the \textbf{points} of $\cC$.
The \textbf{exchange graph} of $\cC$ is a labeled non-oriented graph with vertices set $\cX$
and edges set $I$,
where two vertices $X, Y$ are connected by an edge $i$ if and only if $X\not=Y$ and $r_i(X)=Y$ (and $r_i(Y)=X$).
We display one edge with several labels instead of several edges for simplification.
We say that $\cC$ is \textbf{connected} if its
exchange graph is connected.
For any point $X\in \cX$ we obtain that $A^X$ is not decomposable if $\cC$ is connected by~\cite[Proposition~4.6]{c-Heck09b}.

Recall that there exists a unique category $\cD(\cX, I)$
with objects Ob$\cD(\cX, I)=\cX$ and morphisms
$\Hom(X,Y)=\{(Y,f,X)|f\in \End(\ndZ ^I)\}$
for $X, Y\in \cX$
with the composition defined by
$(Z,g,Y)\circ (Y,f,X)=(Z, gf, X)$
for all $ X, Y, Z\in \cX$,
$f,g\in \textrm{End}(\ndZ ^I)$.

We fix once and for all the notation $(\alpha_i)_{i\in I}$ for the standard basis of $\ndZ^I$.
For all $X\in \cX$ and $i, j\in I$,
let
\begin{equation}\label{eq-si}
s_i^X\in \Aut(\ndZ^I), ~~
s_i^X \alpha _j=\alpha_j-a_{ij}^X \alpha_i.
\end{equation}
Let \textbf{$\cW(\cC)$} be the smallest subcategory of $\cD(\cX, I)$,
where the morphisms are generated by $(r_i(X), s_i^X, X)$,
where $i\in I$, $X\in \cX$.
We will write $s_i^X$ instead of $(r_i(X), s_i^X, X)$, if no confusion is possible.
Notice that $\cW(\cC)$ is a groupoid since all generators $s_i^X$ are invertible.

For all $X\in \cX$,
we say that
$$\rersys{X}=\{w\alpha_i \in \ndZ^I|w\in \cup_{Y\in \cX}\Hom(Y,X)\}$$
is the set of \textbf{real roots} of $X$.
The set $\rersys{X}_{\boldsymbol{+}}=\rersys{X}\cap \ndN_0^I$ is the set of \textbf{positive real roots} of $X$ and $\rersys{X}_{\boldsymbol{-}}=\rersys{X}\cap -\ndN_0^I$ the set of \textbf{negative real roots} of $X$.
The semi-Cartan graph $\cC$ is \textbf{finite} if the set $\rersys{X}$ is finite for all $X\in \cX$.

The semi-Cartan graph $\cC$ is a \textbf{Cartan graph} if the following hold:
\begin{enumerate}
\itemsep=0pt
  \item[(1)] For all $X\in \cX$ the set $\rersys{X}=\rersys{X}_{\boldsymbol{+}}\bigcup \rersys{X}_{\boldsymbol{-}}$.
  \item[(2)] If $l_{mn}^Y:=|\rersys{Y}\cap (\ndN_0 \alpha_m+\ndN_0 \alpha_n)|$ is finite, then $(r_m r_n)^{l_{mn}^Y}(Y)=Y$, where $m, n\in I$, $Y\in \cX$.
\end{enumerate}
In this case,
we say that $\cW(\cC)$ the \textbf{Weyl groupoid of $\cC$} if $\cC$ is a Cartan graph.

We say that
$\cR=\cR(\cC,(\rsys ^{X})_{X\in \cX})$
is a \textbf{root system of type $\cC$}
if for all $X \in \cX $,
the sets $\rsys ^X \subset \ndZ ^I$ are subsets such that:
\begin{enumerate}
\itemsep=0pt
\item[(1)] $\rsys ^X=(\rsys^ X\cap \ndN _0^I)\cup -(\rsys ^X\cap \ndN _0^I)$.
\item[(2)] $\rsys ^X\cap \ndZ \al _i=\{\al _i,-\al _i\}$ for all $i\in I$.
\item[(3)] $s_i^X(\rsys ^X)= \rsys ^{r_i(X)}$ for all $i \in I$.
\item[(4)] $(r_i r_j)^{m_{ij}^X}(X) = X$ for any $i,j \in I$ with $i\not=j$ where
  $m_{ij}^X= |\rsys^X\cap (\ndN _0\al _i + \ndN _0 \al_j)|$ is finite.
\end{enumerate}
\begin{remark}\label{re:irre}
\begin{enumerate}
  \item For any finite Cartan graph $\cC$,
there is a unique root system $\cR=\cR(\cC,(\rersys{X})_{X\in \cX})$ of type $\cC$,
see~\cite[Propositions~2.9,~2.12]{c-Heck09b}.
\item The concept of irreducible root system of type $\cC$ was introduced in~\cite[Definition~4.3]{c-Heck09b}.
If $\cC$ is a finite connected indecomposable Cartan graph,
then the root system $\cR=\cR(\cC,(\rersys{X})_{X\in \cX})$ is an irreducible root system of type $\cC$ by~\cite[Propositions~4.6]{c-Heck09b}.
\end{enumerate}
\end{remark}

\subsection{Cartan graphs for Nichols algebras of diagonal type}\label{se:reCaNI}
Let $\Bbbk G$ be the group algebra over an abelian group $G$.
Let $\ydD$ be the category of
\YD modules over $\Bbbk G$.
We say that $V\in \ydD$ is a \YD
module over $G$.
Write $\ffg$ for the set of $\theta $-tuples of finite-dimensional
irreducible objects in $\ydD$ and $\fiso^G$ for the set of $\theta $-tuples of isomorphism classes
of finite-dimensional irreducible objects in $\ydD$.

We say that a \YD module $V\in \ydD$ is \textbf{of diagonal type} if
$V$ can be decomposed into the direct sum of one-dimensional \YD submodules over $G$.
For any \YD module of diagonal type,
it could be realized as a braided vector space of diagonal type.
Assume that $V$ is of diagonal type.
There exist numbers $q_{ij} \in \Bbbk^*$,
a basis $\{x_i|i \in I\}$
of $V$, and elements $\{g_i|i\in I\}\subset G$ such that
\begin{equation}\label{eq-coaction}
\lcoa (x_i)=g_i\ot x_i, \quad g_i\lact x_j=q_{ij}x_j
\end{equation}
for all $i,j\in I$.
We obtain that the braiding $c(x_i\ot x_j)=q_{ij}x_j\ot x_i$
for all $i,j\in I$.
In this case,
we say that the braiding $c \in \End (V\ot V)$ of $V$ is of diagonal type,
see~\cite[Def.\,5.4]{inp-Andr02}.
The Nichols algebra $\cB(V)$ is said to be~\textbf{of diagonal type}~\cite[Definition~5.8]{inp-Andr02}.
By the proof of~\cite[Proposition~1.3]{HS101},
the matrix $(q_{ij})_{i,j\in I}$ does not depend on the choice of the basis $\{x_i|i\in I\}$,
up to a permutation of $I$.
The matrix $(q_{ij})_{i,j\in I}$ is the~\textbf{braiding matrix of $V$}.
We say that the braiding matrix $(q_{ij})_{i,j\in I}$ is \textbf{indecomposable} if for all $i\not= j$,
there exists a sequence $i_1 = i, i_2, \dots, i_t = j$ of elements of $I$
such that $q_{i_si_{s+1}}q_{i_{s+1}i_s}\not= 1$, $1\leq s\leq t-1$.
Otherwise, we say that the matrix
is decomposable.
Then we can refer to the components of the matrix.

For the remaining part of this section,
let $V$ be any braided vector space of diagonal type with an indecomposable braiding matrix $(q_{ij})_{i, j\in I}$ of a basis $\{x_i|i\in I\}$.
Choose an abelian group $G_0$ and elements $\{g_i|i\in I\}\subset G_0$ such that the assignments in~(\ref{eq-coaction}) define a \YD module structure on $V$.
Then $V=\oplus_{i\in I}\Bbbk x_i$ and each $\Bbbk x_i$ is one-dimensional \YD modules over $G_0$.
Hence $V$ is also a \YD module of diagonal type over $G_0$.

The \textbf{Dynkin diagram} of $V$~\cite[Definition~4]{a-Heck04e}
is a labeled and non-oriented graph $\cD$ satisfying the conditions:
\begin{enumerate}
\item[(1)] For all $i\in I$ the $i$-th vertex is labeled by $q_{ii}$.
\item[(2)] For all $i,j\in I$ there are $n_{ij}$ edges between the $i$-th and $j$-th vertex.
           If $i=j$ or $q_{ij}q_{ji}=1$ then $n_{ij}=0$, otherwise $n_{ij}=1$ and the edge is labeled by $q_{ij}q_{ji}.$
\end{enumerate}
Actually,
we can get the Dynkin diagram of $V$ from the braiding matrix of $V$.

For any tuple $M=(M_j)_{j\in I}\in \ffg$,
write $\cB(M):=\cB(\oplus_{j\in I} M_{j})$ for the Nichols algebra of $M$.
For $x\in \oplus_{j\in I} M_{j}$, $y\in \cB(M)$,
we write $ad_cx(y)=xy-(x_{(-1)}\cdot y)x_{(0)}$ for the adjoint action in the braided category $\ydD$
where $\delta(x)=x_{(-1)}\otimes x_{(0)}$ is the left coaction of $\Bbbk G$.
Notice that
Nichols algebra $\cB(M)$ is $\ndN_0^{\theta}$-graded with
$\deg M_i=\al_i $ for all $i\in I$.
Following the terminology in~\cite{HS10},
we say that the
Nichols algebra $\cB(M)$ is \textbf{decomposable}
if there exists a totally ordered index set $(L,\le)$ and a sequence
$(W_l)_{l\in L}$ of finite-dimensional irreducible $\ndN _0^\theta $-graded objects in
$\ydD $
such that
\begin{equation}\label{eq-decom}
  \cB (M)\simeq
  \bigotimes _{l\in L}\cB (W_l).
\end{equation}
For each decomposition~(\ref{eq-decom}),
we define the set of~\textbf{positive roots} $\rsys^{[M]}_{+}\subset \ndZ^I$ and the set of~\textbf{roots} $\rsys^{[M]}\subset \ndZ^I$ of $[M]$ by
$$\rsys^{[M]}_{+}=\{\deg(W_l)|\, l\in L\}, \quad
\rsys^{[M]}=\rsys^{[M]}_{+}\cup-\rsys^{[M]}_{+}.$$
By~\cite[Theorem~4.5]{HS10} we obtain that
the set of
roots $\rsys^{[M]}$ of $[M]$ does not depend on the choice of the decomposition.
\begin{remark}\label{rem-decom}
 If $\dim M_i=1$ for all $i\in I$,
 then the Nichols algebra $\cB(M)$ is decomposable based on the theorem of V.~Kharchenko~\cite[Theorem~2]{a-Khar99}.
\end{remark}

Let $i\in I$.
Following the terminology in~\cite[Definition~6.4]{HS10} we say that $M$ is \textbf{$i$-finite},
if there is $m\in \ndN$ such that $(\ad_{c} M_i)^m (M_j)=0$ for any $j\in I\setminus \{i\}$.
For all $i \in I$ and any $j\in I \setminus \{i\}$ let $a_{ij}^{M}=-\infty$ if $M$ is not $i$-finite and let
$a_{ij}^M=-\mathrm{max}\{m\in \ndN_0 \,|\,(\ad_c  M_i)^m(M_j)\not=0 \}$ otherwise.
Let also $a_{ii}^M=2$ for all $i\in I$.
The matrix $(a_{ij}^M)_{i,j\in I}$ is called the~\textbf{Cartan matrix} of $M$,
denoted by $A^M$.
For all $i \in I$ let $R_i(M)=M$ if $M$ is not $i$-finite and
let $R_i(M)=({R_i(M)}_j)_{j\in I}$, where
\begin{align}\label{eq-ri}
{R_i(M)}_j=&
  \begin{cases} {M_i}^* & \text{if $j=i$,}\\
(\ad_{c}M_i)^{-a_{ij}^M}(M_j) & \text{if $j\not=i$.}
  \end{cases}
\end{align}
The tuple $R_i(M)\in \ffg$ since $R_i(M)_j$ is irreducible by~\cite[Theorem~7.2(3)]{HS10}.
Then $R_i$ is a reflection map from $\ffg$ to $\ffg$.
We say that $R_i$ is the i-th reflection.

For any tuple $M\in\ffg$,
we define
$$\fiso^G(M)=\{[R_{i_1} \cdots R_{i_n}(M)]\in \fiso^G |\,n\in \ndN_0, i_1,\dots, i_n\in I\}$$
$$\ffg(M)=\{R_{i_1} \cdots R_{i_n}(M)\in \ffg|\, n\in \ndN_0, i_1,\dots, i_n\in I\}.$$

We say that $M$ \textbf{admits all reflections} if $N$ is $i$-finite for all $N\in \ffg(M)$.

Recall that $V$ is a braided vector space of diagonal type with the indecomposable braiding $c(x_i\ot x_j)=q_{ij}x_j\ot x_i$
for all $i,j\in I$.
For the remaining part of this section,
we consider the tuple $M:=(\Bbbk x_1, \Bbbk x_2, \dots, \Bbbk x_{\theta})$ of one-dimensional \YD modules over $G_0$.
We say that the braiding matrix of $V$ is the \textbf{braiding matrix of $M$} and that the Dynkin diagram of $V$ is the \textbf{Dynkin diagram of $M$}.
Since the braiding matrix $(q_{ij})_{i, j\in I}$ is indecomposable,
the tuple $M$ is braid-indecomposable,
see~\cite[Definition~2.1]{VH-14}.

Let $i\in I$.
Assume that $M$ is $i$-finite.
Then $y_j$ is a basis of ${R_i(M)}_j$ defined by~(\ref{eq-ri}),
where
\begin{align}\label{eq-yi}
y_j:=&
 \begin{cases}
 y_i                               & \text{if $j=i$,}\\
 (\ad_{c} x_i)^{-a_{ij}^M}(x_j) & \text{if $j\not=i$,}
 \end{cases}
\end{align}
where $y_i\in M_i^*\setminus \{0\}$.

The following Lemma~\ref{lem:aij} gives us the way to check whether $M$ is $i$- finite, for $i\in I$,
It can be obtained from~\cite[Lemma~3.7]{inp-AndrSchn02},
see~\cite[Lemma~3.2]{WH-14} for details.
\begin{lemma}\label{lem:aij}
Let $i\in I$.
Then the tuple $M$ is $i$-finite if and only if for any $j\in I\setminus\{i\}$
there is an integer $m\in \ndN_0$ satisfying $(m+1)_{q_{ii}}(q_{ii}^mq_{ij}q_{ji}-1)=0$.
In this case,
\begin{align}\label{eq-aij}
  a_{ij}^M=-\mathrm{min}\{m\in \ndN_0 \,|\,(m+1)_{q_{ii}}(q_{ii}^mq_{ij}q_{ji}-1)=0\}.
\end{align}
\end{lemma}
Further,
we get the labels of the Dynkin diagram of $R_i(M)$ from Lemma~\ref{lem:qij'}.
The method was described in~\cite[Example~1]{a-Heck04e} and~\cite[Lemma~3.4]{WH-14}.
\begin{lemma}\label{lem:qij'}
Let $i\in I$.
Assume that $M$ is $i$-finite and $(a_{ij})_{j\in I}:=(a_{ij}^M)_{j\in I}$ defined by~(\ref{eq-aij}).
Let $(l_{jk})_{j,k\in I}$ be the braiding matrix of $R_i(M)$ with respect to the basis $(y_j)_{j\in I}$ defined by~(\ref{eq-yi}).
Write $q_{ij}':=q_{ij}q_{ji}$ and $l_{ij}':=l_{ij}l_{ji}$ for all $i,j\in I$.
Then%
\begin{align*}
l_{jj}=
 \begin{cases}
 q_{ii} &\text{if}~j=i,\\
 q_{jj}{(q_{ij}')^{-a_{ij}}} &\text{if}~ j\not=i, q_{ii}=1,\\
 q_{jj}  &\text{if}~j\not=i, q_{ij}'=q_{ii}^{a_{ij}},\\
 q_{ii}q_{jj}{(q_{ij}')^{-a_{ij}}} &\text{if}~ j\not=i, q_{ii}\in G'_{1-a_{ij}},
 \end{cases}
 &
 ~l_{ij}'=
 \begin{cases}
 (q_{ij}')^{-1} & \text{if}~ j\not=i, q_{ii}=1,\\
 q_{ij}' &\text{if} ~j\not=i, q_{ij}'=q_{ii}^{a_{ij}},\\
 q_{ii}^2(q_{ij}')^{-1} &\text{if}~ j\not=i, q_{ii}\in G'_{1-a_{ij}},
 \end{cases}
\end{align*}
 and
 \begin{align*}
l_{jk}'=
 \begin{cases}
 q_{jk}'(q_{ij}')^{-a_{ik}}(q_{ik}')^{-a_{ij}} &  \text{if} ~q_{ii}=1,\\
 q_{jk}'  & \text{if}~q_{ir}'=q_{ii}^{a_{ir}}, r\in \{j, k\},\\
 q_{jk}'(q_{ik}'q_{ii}^{-1})^{-a_{ij}}& \text{if} ~q_{ii}\in G'_{1-a_{ik}}, q_{ii}^{a_{ij}}=q_{ij}', \\
  q_{jk}'q_{ii}^{2}(q_{ij}'q_{ik}')^{-a_{ij}} & \text{if} ~q_{ii}\in G'_{1-a_{ik}}, q_{ii}\in G'_{1-a_{ij}}.
  \end{cases}
 \end{align*}
 for any $j, k\not=i$, $j\not=k$.
\end{lemma}
Hence we can check whether $M$ admits all reflections by Lemmas~\ref{lem:aij} and~\ref{lem:qij'}.

Assume that the tuple $M$ admits all reflections.
The following Proposition implies that we can attach to $M$ a finite connected Cartan graph $\cC(M)$.
\begin{prop}\label{prop:small}
  Assume that the tuple $M$ admits all reflections and that $\cW(M):=\cW(\cC(M))$ is defined.
  Let the tuple $M$ be braid-indecomposable and $\cW(M)$ is finite.
  For all $X\in \fiso^{\Bbbk G_0}(M)$
  let $$[X]_{\theta}^s=\{Y\in \fiso^{\Bbbk G_0}(M)| \,\text{Y and X have the same Dynkin diagram}\}$$
  and $\cY_{\theta}^s(M)=\{[X]^s_{\theta} |\, X\in \fiso^{\Bbbk G_0}(M)\}.$
  Let $A^{[X]^s_{\theta}}=A^X$ for all $X\in \fiso^{\Bbbk G_0}(M)$ and
  $$t: I\times \cY_{\theta}^s(M)\rightarrow \cY_{\theta}^s(M), \quad
  (i, [X]^s_{\theta})\mapsto [R_i(X)]^s_{\theta}.$$
  Then the tuple
  $$\cC(M)=\{I, \cY_{\theta}^s(M), t, (A^Y)_{Y\in \cY_{\theta}^s(M)}\}$$
  is a finite connected indecomposable Cartan graph.
\end{prop}
\begin{proof}
 Since the tuple $M$ admits all reflections,
 we obtain that $\cC(M)$ is a well-defined connected semi-Cartan graph by~\cite[Proposition~3.6]{WH-14}.
 By~\cite[Theorem~2.3 and Corollary~2.4]{HS101} $\cC(M)$ is a finite Cartan graph since $\cW(\cC(M))$ is finite.
 The braid-indecomposability of $M$ implies that the Cartan matrix $A^{M}$ of $M$ is indecomposable.
 Hence $\cC(M)$ is indecomposable.
\end{proof}
\section{Cartan graphs of rank three}\label{se:Rank3Cartan}
In this section,
we obtain Proposition~\ref{prop:A3},
which illustrates the properties of finite connected indecomposable Cartan graphs of rank three.
Proposition~\ref{prop:A3} will be used for our classification in the next section.
Our main referee is~\cite{c-Heck12a}.

Let $\cC = \cC (I, \cX, r, (A^X)_{X \in \cX})$ be a semi-Cartan graph where $I=\{1,2,3\}$.
Let $X\in\cX$.
Since the points of $\cC$ could have many different neighborhoods,
we define the following good neighborhoods in order to cover all the finite connected indecomposable Cartan graphs in such way that
at least one point of $\cC$ has one of the good neighborhoods.
\begin{defn}\label{defA3}
  We say that $X$ has a \textbf{good $A_3$ neighborhood} if there exists integer sequence $(a,b,c,d)$ such that there is a permutation of $I$ with respect to which
  $A^X=\begin{pmatrix}2&-1&0\\  -1&2&-1\\  0&-1&2  \end{pmatrix}$,
  $A^{r_1(X)}=\begin{pmatrix}2&-1&0\\-1&2&-a\\0&-1&2\end{pmatrix}$,
  $A^{r_2(X)}=\begin{pmatrix}2&-1&-b\\-1&2&-1\\-c&-1&2\end{pmatrix}$,
  $A^{r_3(X)}=\begin{pmatrix}2&-1&0\\-d&2&-1\\0&-1&2\end{pmatrix}$,
  and one of the following holds.
  \begin{enumerate}
    \item[(1)] $(a,b,c,d)\in \{(1,0,0,1), (1,1,1,1), (1,1,1,2), (1,1,2,3), (2,1,1,2), (2,1,2,2)\}$.
    \item[(2)] $(a,b,c,d)=(2,1,2,3)$, $-a_{21}^{r_3r_1(X)}=3$.
  \end{enumerate}
\end{defn}
\begin{defn}\label{defB3}
  We say that $X$ has a \textbf{good $B_3$ neighborhood} if there is a permutation of $I$ with respect to which\\
  $A^X=A^{r_1(X)}=A^{r_2(X)}=\begin{pmatrix} 2&-1&0\\ -1&2&-1\\ 0&-2&2 \end{pmatrix}$, $A^{r_3(X)}=\begin{pmatrix} 2&-1&0\\ -a&2&-1\\ 0&-2&2 \end{pmatrix}$, and $-a^{r_1r_3(X)}_{23}\in \{1,2\}$, where $a\in \{1,2\}$.
\end{defn}
\begin{defn}\label{defC3}
  We say that $X$ has a \textbf{good $C_3$ neighborhood} if there is a permutation of $I$ with respect to which\\
  $A^X=A^{r_1(X)}=A^{r_2(X)}=\begin{pmatrix} 2&-1&0\\ -1&2&-2\\ 0&-1&2 \end{pmatrix}$
  and $A^{r_3(X)}=\begin{pmatrix}2&-1&0\\-a&2&-2\\0&-1&2\end{pmatrix}$,
  where $a\in \{1, 2\}$.
\end{defn}
Based on Theorem~4.1 in~\cite{c-Heck12a},
we get the following result by computer calculations algorithms.
\begin{prop}\label{prop:A3}
 Let $\cC = \cC (I,\cX,r,(A^X)_{X \in \cX})$ be a rank three finite connected indecomposable Cartan graph.
 Then there exists at least one point $Y\in \cX$ satisfying one of the following:
 \begin{enumerate}
   \item[(a)] The point $Y$ has a good $A_3$ neighborhood.
   \item[(b)] The point $Y$ has a good $B_3$ neighborhood.
   \item[(c)] The point $Y$ has a good $C_3$ neighborhood.
 \end{enumerate}
\end{prop}

\begin{remark}\label{rem-good}
 There are finite connected indecomposable Cartan graphs of rank three which contain precisely one case of the good neighborhoods.
\end{remark}
\begin{proof}
Since $\cC$ is a finite connected indecomposable Cartan graph,
$\cR=\cR(\cC,(\rersys{X})_{X\in \cX})$
is a finite irreducible root system of type $\cC$,
see Remark~\ref{re:irre}.
By~\cite[Theorem~4.1]{c-Heck12a} there exists a point $X\in \cX$ satisfying that the set $\rersys{X}_{\boldsymbol{+}}$ is in the list of~\cite[Appendix A]{c-Heck12a} up to a permutation of $I$.
There are precisely 55 such possible sets for the rank 3.
We consider each set in the list.
For each point $Y$ in the list,
we calculate all neighbors $\{r_i(Y)$, $i\in I$\} of $Y$.
Since the reflection $s_i^Y$ maps $\rersys{Y}_{\boldsymbol{+}}\setminus \{\alpha_i\}$ bijectively to $\rersys{r_i(Y)}_{\boldsymbol{+}}\setminus \{\alpha_i\}$ for any $i\in I$ ,
the Cartan matrices of all neighbors of $Y$ can be obtained from $\rersys{Y}_{\boldsymbol{+}}$, $Y\in \cX$ by~\cite[Corollary~2.9]{c-Heck12a}.
If $Y$ has a good $A_3$, $B_3$ or $C_3$ neighborhood,
then the proof is done.
Otherwise repeat the previous step to the neighbours of $Y$.
Since $\cX$ is finite,
this algorithm terminates.
The elementary calculations are done by computer algorithms and they are skipped here.
\end{proof}

\section{The classification result}\label{se:clasi}
In this section we explain the classification of rank three Nichols algebras of diagonal type with a finite root system over fields of positive characteristic.
We formulate the main result in Theorem~\ref{Theo:clasi} and Corollary~\ref{coro-cla}.
All Dynkin diagrams of rank three braided vector spaces with a finite root system are listed in Tables~(\ref{tab.1})-(\ref{tab.3}).
\begin{theorem}\label{Theo:clasi}
Let $\Bbbk$ is a field of characteristic $p>0$.
Let $V$ be a braided vector space of diagonal type with $\dim_{\Bbbk}V=3$.
Let $I=\{1,2,3\}$.
Let $\{x_k|k\in I\}$ be a basis of $V$ such that there exists an indecomposable braiding matrix
$(q_{ij})_{i,j\in I}\in (\Bbbk^*)^{3\times 3}$ with respect to $\{x_k|k\in I\}$.
Let the tuple $M:=(\Bbbk x_i)_{i\in I}$.
Then the following are equivalent.
  \begin{enumerate}
  \item[(1)] $M$ admits all refections and $\cW(\cC(M))$ is finite.
  \item[(2)] The Dynkin diagram $\cD$ of $V$ appears in the Tables~\ref{tab.1}, \ref{tab.2}, and \ref{tab.3}, if $p=2$, $p=3$, and $p>3$, respectively.
  \item[(3)] The Nichols algebra $\cB(V)$ is decomposable and the set of roots ${\roots}^{[M]}$ is finite.
 \end{enumerate}
In this case,
the row of Tables~(\ref{tab.1})-(\ref{tab.3}) containing $\cD$
consists precisely of the Dynkin diagrams of the points of $\cC(M)$.
The corresponding row of Table~\ref{tab.4} contains the exchange graph of $\cC(M)$.
\begin{remark}
\begin{enumerate}
\item[(1)] We present the i-vertex of $\cD$ with a bullet point $\bullet$ to describe the i-th reflection map $R_i$ of the Dynkin diagram $\cD$.
 \item[(2)] In order to give the exchange graph of $\cC(M)$ in Theorem~\ref{Theo:clasi},
  we use the following notations in Tables~(\ref{tab.1})-(\ref{tab.3}).
  For each row $n$ of Tables~(\ref{tab.1})-(\ref{tab.3}),
  we enumerate $l$-th Dynkin diagram with $\cD_{nl}$ column by column for all $l\geq1$.
  For each Dynkin diagram,
  we enumerate the vertices from left to right with $1, 2$ and $3$, respectively.
  For $i,j,k\in I$,
  we also write $\tau_{ijk} \cD_{nl}$ for the graph $\cD_{nl}$
  where the three vertices of $\cD_{nl}$ change to $i$, $j$, $k$, respectively.
   \end{enumerate}
\end{remark}
\end{theorem}
\begin{proof}
  The equivalence between $(1)$ and $(3)$ is clear, see~\cite[Theorem~5.1]{WH-14}.

  (2)$\Rightarrow$ (1).
  Assume that the Dynkin diagram $\cD$ appears in row $r$ of any of Tables~(\ref{tab.1})-(\ref{tab.3}).
  From Lemmas~\ref{lem:aij} and~\ref{lem:qij'} we determine that $M$ admits all reflections.
  Hence the Cartan graph $\cC(M)$ can be defined by Proposition~\ref{prop:small}.
  We obtain that $\cC(M)$ is the same with the Cartan graph obtained from the Dykin diagrams appearing in row $s$ of Table~2 in~\cite{a-Heck04aa},
  where $s$ appears in the third column of row $r$ of Table~\ref{tab.4}.
  The detailed calculations are skipped at this point here.
  Moreover,
  the Weyl groupoids of the above Cartan graphs from~\cite{a-Heck04aa} are finite,
  see~\cite[Table~2]{c-Heck12a}.
  Hence $\cW(\cC(M))$ is finite.

  (1)$\Rightarrow$ (2).
  Since $M$ admits all reflections,
  we can define a semi-Cartan graph $\cC(M)$ by Proposition~\ref{prop:small}.
  Hence $\cW(M)=\cW(\cC(M))$ is defined.
  Further,
  $\cC(M)$ is a finite connected indecomposable Cartan graph since $\cW(\cC(M))$ is finite and the matrix $(q_{ij})_{i,j\in I}$ is indecomposable by Proposition~\ref{prop:small}.
  From the implication (2)$\Rightarrow$ (1) it is enough to prove that
  there exists at least one point in $\cC(M)$ contained in Tables~(\ref{tab.1})-(\ref{tab.3}).

  Set $X=[M]^s_3$, $q_1=q_{12}q_{21}$, $q_2=q_{23}q_{32}$, and $q_3=q_{31}q_{13}$.
  We write the Cartan matrix $A^X:=(a_{ij})_{i,j\in I}$.
  Since $\cC(M)$ is a finite Cartan graph,
  by Proposition~\ref{prop:A3} we are free to assume that Cartan matrix $A^X$
  is one of the following cases: (a), (b) and (c).

  Case (a). Assume that $X$ has a good $A_3$ neighborhood.
  Let $$(a,b,c,d):=(-a_{23}^{r_1(X)}, -a_{13}^{r_2(X)}, -a_{31}^{r_2(X)}, -a_{21}^{r_3(X)})$$ be the sequence in Definition~\ref{defA3}.
  Applying Lemma~\ref{lem:aij} to $a_{13}=0$, $a_{12}=-1$, and $a_{23}=-1$,
  we get $q_3=1$, $q_1\not=1$ and $q_2\not=1$, respectively.
  Since $A^X$ is of $A_3$ type,
  we distinguish subcases $a_1$, $a_2$, $a_3$, $a_4$, $a_5$, $a_6$, $a_7$, and $a_8$ by Lemma~\ref{lem:aij}.

  Case $a_1$.
  Consider the case $(2)_{q_{11}}=(2)_{q_{22}}=(2)_{q_{33}}=0$.
  Set $q_1:=q$ and $q_2:=r$.
  Then $q_{11}=q_{22}=q_{33}=-1$.
  If $qr=1$ and $q\not=-1$,
  then $\cD=\cD_{82}$.
  If $qr=1$ and $q=-1$,
  then $p\not=2$ and $\cD=\cD_{11}$, where $q=-1$, and $p\not=2$.
  If $qr\not=1$,
  then by Lemma~\ref{lem:qij'} we get the reflection $r_2(X)$ of $X$
  \setlength{\unitlength}{1mm}
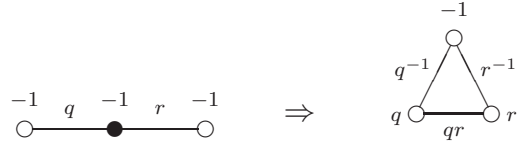
\begin{gather*}
\Dchainthree{m}{$-1$}{$q$}{$-1$}{$r$}{$-1$} \quad \quad \Rightarrow \quad \quad
\rule{0pt}{12\unitlength}
\Dtriangle{}{$q$}{$-1$}{$r$}{$q^{-1}$}{$r^{-1}$}{$qr$}\quad
\end{gather*}
 with $q\not=1$, $r\not=1$, and $qr\not=1$.
 Since $X$ has a good $A_3$ neighborhood,
 we obtain that $a_{13}^{r_2(X)}=-1$ and $a_{31}^{r_2(X)}\in \{-1,-2\}$.
 Then we distinguish subcases $a_{1a}$ and $a_{1b}$.

  Case $a_{1a}$.
  Consider the case $a_{13}^{r_2(X)}=a_{31}^{r_2(X)}=-1$.
  We obtain that $(2)_q(q^2r-1)=0$ and $(2)_r(qr^2-1)=0$ by Lemma~\ref{lem:aij}.
  If $q=-1$, then $p\not=2$ and $r\not=-1$ by $qr\not=1$.
  If $q=-1$ and $qr^2-1=0$, then $p\not=2$, $r\in G'_4$ and hence $\cD=\tau_{321}\cD_{62}$, where $q\in G'_4$ and $p\not=2$.
  If $q^2r=1$ and $r=-1$, then $p\not=2$ and $\cD=\cD_{62}$, where $q\in G'_4$ and $p\not=2$.
  If $q^2r=qr^2=1$, then $q=r\in G'_3$, $p\not=3$ and hence $\cD=\cD_{15,2}$, $p\not=3$.

  Case $a_{1b}$.
  Consider the case $a_{13}^{r_2(X)}=-1$ and $a_{31}^{r_2(X)}=-2$.
  Then $(2)_q(q^2r-1)=0$, $(3)_r(qr^3-1)=0$ and $(2)_r(qr^2-1)\not=0$.
  Hence $r\not=-1$ and $qr^2\not=1$.
  If $q=-1$, $r\in G'_3$ and $p\not=3$,
  then $p\not=2$ and $\cD=\cD_{17,1}$, $p\not=2,3$.
  If $q=-1$ and $qr^3=1$,
  then $r\in G'_6$, $p\not=2,3$, and hence $\cD=\tau_{321}\cD_{72}$, where $q\in G'_6$ and $p\not=2, 3$.
  If $q^2r=r^3=1$,
  then $-q=r\in G'_3$, $p\not=2,3$ by $qr^2\not=1$.
  Hence $\cD=\tau_{321}\cD_{17,3}$, $p\not=2, 3$.
  If $q^2r=qr^3=1$,
  then $r^2=q$, where $r\in G'_5$, $p\not=5$.
  Hence by Lemmas~\ref{lem:aij} and~\ref{lem:qij'} we get the sequence $(a,b,c,d)=(2,1,2,3)$.
  Further $a_{21}^{r_3r_1(X)}=-4$,
  which is a contradiction.

  Case $a_2$.
  Consider the case $(2)_{q_{22}}=(2)_{q_{33}}=0$, $q_{11}q_1-1=0$, and $(2)_{q_{11}}\not=0$.
  Set $q_{11}:=q$ and $q_2:=r$.
  Then we get $q\not=-1$ for $p\not=2$ and $q\not=1$ for $p=2$.
  Hence we obtain $q_{22}=q_{33}=-1$, $qq_1=1$, $q\not=-1$ for $p\not=2$ or $q_{22}=q_{33}=1$, $qq_1=1$, $q\not=1$ for $p=2$.
  We distinguish two subcases $a_{2a}$ and $a_{2b}$.

  Case $a_{2a}$
  Consider the case $q_{22}=q_{33}=-1$, $qq_1=1$, $q\not=-1$, and $p\not=2$.
  If $r=q\not=-1$,
  then $\cD=\tau_{321}\cD_{42}$, $p\not=2$.
  If $r=-1$, $q\not=r$, and $q\notin G_4'$,
  then $\cD=\cD_{91}$, where $r=-1$, $q\notin G_4'$ and $p\not=2$.
  If $r=-1$, $q\not=r$, and $q\in G_4'$,
  then $\cD=\cD_{10,2}$, where $q\in G_4'$ and $p\not=2$.
  Now we consider the case $r\not=-1$ and $q\not=r$.
  By Lemma~\ref{lem:qij'} we get the following reflection $r_2(X)$ of $X$
  \setlength{\unitlength}{1mm}
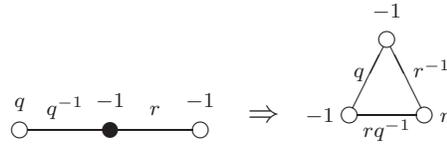
\begin{gather*}\label{dia2}
\Dchainthree{m}{$q$}{$q^{-1}$}{$-1$}{$r$}{$-1$} \quad \Rightarrow \quad
\rule{0pt}{10\unitlength}
\Dtriangle{}{$-1$}{$-1$}{$r$}{$q$}{$r^{-1}$}{$r q^{-1}$}\quad
\end{gather*}
  with $q\not=1$, $r\notin \{-1,1, q\}$, and $p\not=2$.
  Since $X$ has a good $A_3$ neighborhood,
  we get $a_{31}^{r_2(X)}\in \{-1, -2\}$.
  Then we get either $(2)_r(r^2q^{-1}-1)=0$ or $(3)_r(r^3q^{-1}-1)=0$, $(2)_r(r^2q^{-1}-1)\not=0$.
  Hence we can distinguish three subcases $a_{2a1}$, $a_{2a2}$ and $a_{2a3}$.

  Case $a_{2a1}$.
  Consider the case $(2)_r(r^2q^{-1}-1)=0$.
  Since $r\not=-1$ and $p\not=2$,
  we obtain that $r^2q^{-1}=1$ and hence $\cD=\tau_{321}\cD_{62}$, $p\not=2$.

  Case $a_{2a2}$.
  Consider the case $r^3q^{-1}=1$ and $(2)_r(r^2q^{-1}-1)\not=0$.
  Then we get $q=r^3\not=1$ and $r\notin \{1, -1\}$.
  Hence $\cD=\tau_{321}\cD_{72}$, $p\not=2$.

  Case $a_{2a3}$.
  Consider the case $r\in G_3'$, $(2)_r(r^2q^{-1}-1)\not=0$, and $p\not=3$.
  Then $qr\not=1$.
  By Lemma~\ref{lem:qij'} we get the following reflection $r_3(X)$  of $X$
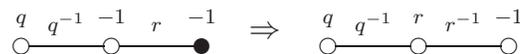
  \begin{gather*}
\rule{0pt}{6\unitlength}
\Dchainthree{r}{$q$}{$q^{-1}$}{$-1$}{$r$}{$-1$} \quad \Rightarrow \quad
\Dchainthree{}{$q$}{$q^{-1}$}{$r$}{$r^{-1}$}{$-1$}\quad \quad
\end{gather*}
  with $q\not=1$, $r\notin \{-1,1, q, q^{-1}\}$.
  By Lemmas~\ref{lem:aij} and~\ref{lem:qij'} we obtain that $(a,b,c,d)=(1,1,2,2)$, which is a contradiction.

  Case $a_{2b}$.
  Consider the case $q_{22}=q_{33}=1$, $qq_1=1$, $q\not=1$ for $p=2$.
  If $r=q$,
  then $\cD=\tau_{321}\cD_{42}$, $p=2$.
  If $r\not=q$,
  then we get $a_{31}^{r_2(X)}\in \{-1, -2\}$.
  Hence we distinguish subcases $a_{2b1}$ and $a_{2b2}$.

  Case $a_{2b1}$.
  If $a_{31}^{r_2(X)}=-1$, then $r^2q^{-1}=1$ and hence $\cD=\tau_{321}\cD_{62}$, $p=2$.

  Case $a_{2b2}$.
  Consider the case $a_{31}^{r_2(X)}=-2$.
  If $r\in G_3'$,
  then we get $(a,b,c,d)=(1,1,2,2)$, which is a contradiction.
  If $r^3q^{-1}=1$,
  then $q=r^3\not=1$ and hence $\cD=\tau_{321}\cD_{72}$, $p=2$.

   Case $a_3$.
  Consider the case $(2)_{q_{11}}=(2)_{q_{22}}=0$, $q_{33}q_2-1=0$, and $(2)_{q_{33}}\not=0$.
  Let $q_{33}:=q$ and $q_1:=r$.
  Then by $(2)_{q}\not=0$ we get $q\not=-1$ for $p\not=2$ and $q\not=1$ for $p=2$.
  Hence we obtain $q_{11}=q_{22}=-1$, $qq_2=1$, $q\not=-1$ for $p\not=2$ or $q_{11}=q_{22}=1$, $qq_2=1$, $q\not=1$ for $p=2$.
  We distinguish two subcases $a_{3a}$ and $a_{3b}$.

  Case $a_{3a}$
  Consider the case $q_{11}=q_{22}=-1$, $qq_2=1$, $q\not=-1$, and $p\not=2$.
  If $r=q\not=-1$,
  then $\cD=\cD_{42}$, $p\not=2$.
  If $r=-1$, $q\not=r$, and $q\notin G_4'$,
  then $\cD=\tau_{321}\cD_{91}$, where $r=-1$, $q\notin G_4'$ and $p\not=2$.
  If $r=-1$, $q\not=r$, and $q\in G_4'$,
  then $\cD=\tau_{321}\cD_{10,2}$, where $q\in G_4'$ and $p\not=2$.
  Now we consider the case $r\not=-1$ and $q\not=r$.
  By Lemma~\ref{lem:qij'} we get the following reflection $r_2(X)$ of $X$
  \setlength{\unitlength}{1mm}
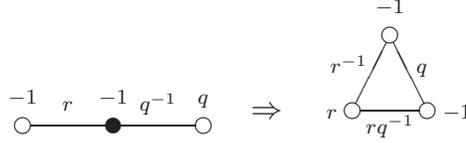
\begin{gather*}\label{dia3}
\Dchainthree{m}{$-1$}{$r$}{$-1$}{$q^{-1}$}{$q$} \quad \Rightarrow \quad
\rule{0pt}{10\unitlength}
\Dtriangle{}{$r$}{$-1$}{$-1$}{$r^{-1}$}{$q$}{$r q^{-1}$}\quad
\end{gather*}
  with $q\not=1$, $r\notin \{-1,1, q\}$, and $p\not=2$.
  Since $X$ has a good $A_3$ neighborhood and $A^{r_3(X)}=A^X$,
  we get $a_{13}^{r_2(X)}=-1$ and hence $r^2q^{-1}=1$.
  Then $\cD=\cD_{62}$, $p\not=2$.

  Case $a_{3b}$.
  Consider the case $q_{11}=q_{22}=1$, $qq_2=1$, $q\not=1$ for $p=2$.
  If $r=q$,
  then $\cD=\cD_{42}$, $p=2$.
  If $r\not=q$,
  then we get $a_{31}^{r_2(X)}=-1$.
  Hence, similarly to the argument in Case $a_{3a}$,  $r^2q^{-1}=1$ and $\cD=\cD_{62}$, $p=2$.

   Case $a_4$.
  Consider the case $(2)_{q_{11}}=(2)_{q_{33}}=0$, $q_{22}q_1-1=q_{22}q_2-1=0$, and $(2)_{q_{22}}\not=0$.
  Let $q_{22}:=q$.
  Then by $(2)_{q}\not=0$ we get $q\not=-1$ for $p\not=2$ and $q\not=1$ for $p=2$.
  Hence we obtain $q_{11}=q_{33}=-1$, $qq_1=qq_2=1$, $q\not=-1$ for $p\not=2$ or $q_{11}=q_{33}=1$, $qq_1=qq_2=1$, $q\not=1$ for $p=2$.
  Hence we get $\cD=\cD_{83}$.

  Case $a_5$.
  Consider the case $(2)_{q_{11}}=0$, $q_{22}q_1-1=q_{22}q_2-1=q_{33}q_2-1=0$, $(2)_{q_{22}}\not=0$, and $(2)_{q_{33}}\not=0$.
  Set $q_{33}:=q$.
  If $q_{11}=1$ and $p=2$, then $\cD=\cD_{41}$, $p=2$.
  If $q_{11}=-1$ and $p\not=2$, then $q\not=-1$ and $\cD=\cD_{41}$, $p\not=2$.

  Case $a_6$.
  Consider the case $(2)_{q_{33}}=0$, $q_{11}q_1-1=q_{22}q_1-1=q_{22}q_2-1=0$, $(2)_{q_{11}}\not=0$, and $(2)_{q_{22}}\not=0$.
  Set $q_{11}:=q$.
  Then $q\not=-1$.
  If $q_{33}=1$ and $p=2$, then $\cD=\tau_{321}\cD_{41}$, $p=2$.
  If $q_{33}=-1$ and $p\not=2$, then $q\not=-1$ and $\cD=\tau_{321}\cD_{41}$, $p\not=2$.

  Case $a_7$.
  Consider the case $(2)_{q_{22}}=0$, $q_{11}q_1=q_{33}q_2=1$, $(2)_{q_{11}}\not=0$, and $(2)_{q_{33}}\not=0$.
  Set $q_{11}:=q$ and $q_{33}:=r$.
  Then $q_{22}=-1$ and $qq_1=rq_2=1$.
  If $qr=1$,
  then $q\not=-1$ and $\cD=\cD_{81}$.
  If $qr\not=1$, $q=r\in G'_3$, then $p\not=3$ and hence $\cD=\cD_{11,1}$, $p\not=3$.
  If $qr\not=1$ and $q=r\notin G'_3$, then $\cD=\cD_{10,1}$.
  If $qr\not=1$, $q\not=r$, $qr^2\not=1$, and $rq^2\not=1$, then $\cD=\cD_{91}$.
  If $qr\not=1$, $q\not=r$, and $rq^2=1$, then $q\notin G_3'$ and $\cD=\cD_{10,2}$.
  If $qr\not=1$, $q\not=r$, and $qr^2=1$, then $r\notin G_3'$, $r^2\not=1$, and $\cD=\tau_{321}\cD_{10,2}$.


  Case $a_8$.
  Consider the case $q_{11}q_1=q_{22}q_1=q_{22}q_2=q_{33}q_2=1$, $(2)_{q_{11}}\not=0$, $(2)_{q_{22}}\not=0$, and $(2)_{q_{33}}\not=0$.
  Then $\cD=\cD_{11}$, $q\not=-1$.

  Case (b).
  Assume that $X$ has a good $B_3$ neighborhood.
  Since $A^X$ is of $B_3$ type,
  we obtain that $q_3=1$, $q_1\not=1$ and $q_2\not=1$ by Lemma~\ref{lem:aij}.
  By Lemma~\ref{lem:qij'}
  we get $q_1q_2=1$ if $(2)_{q_{22}}=0$,
  since $a_{13}^{r_2(X)}=0$.
  Further,
  we distinguish subcases $b_1$, $b_2$, $b_3$, $b_4$, $b_5$, $b_6$, $b_7$, and $b_8$.

  Case $b_{1}$.
  Consider the case $(2)_{q_{11}}=(2)_{q_{22}}=0$, $q_{33}^2q_2-1=0$, and $(2)_{q_{33}}(q_{33}q_2-1)\not=0$.
  Then $q_{11}=q_{22}=-1$.
  Let $q_{33}:=q$.
  Then $q_2=q^{-2}$ and hence $q_1=q^2$ by $q_1q_2=1$.
  If $q^2\not=-1$,
   then $\cD=\cD_{52}$.
   If $q^2=-1$,
   then $p\not=2$ and $\cD=\cD_{21}$, where $q\in G'_4$, $p\not=2$.

  Case $b_{2}$.
  Consider the case $(2)_{q_{11}}=(2)_{q_{22}}=(3)_{q_{33}}=0$,  $q_{33}^2q_2-1\ne0$, and $(2)_{q_{33}}(q_{33}q_2-1)\not=0$.
  Let $q_1:=q$.
  Then $q_2=q^{-1}$ and $q_{33}\not=q$.
  We distinguish subcases $b_{2a}$, $b_{2b}$ and $b_{2c}$.

  case $b_{2a}$.
  Consider the case $q_{11}=q_{22}=-1$, $q_{33}:=\zeta\in G_3'$.
  Then we get $p\not=2, 3$.
  Then $\zeta q\ne 1$.
  By Lemma~\ref{lem:qij'} we get the reflection $r_3(X)$ of $X$
   \begin{gather*}
\rule{0pt}{6\unitlength}
\Dchainthree{r}{$-1$}{$q$}{$-1$}{$q^{-1}$}{$\zeta$} \quad \Rightarrow \quad
\Dchainthree{}{$-1$}{$q$}{-$\zeta q^{-2}$}{$~q{\zeta}^{-1}$}{$\zeta$}\quad \quad
\end{gather*}
   with $\zeta\in G'_3$, $q\not=\zeta,\zeta^{-1}$, and $p\not=2, 3$.
   Since $a_{23}^{r_3(X)}=-1$ by Definition~\ref{defB3},
   we get $((-\zeta q^{-2})q\zeta^{-1}-1)(\zeta q^{-2}-1)=0$ and hence $q=-\zeta^{-1}$ or $q=-1$.
   If $q=-\zeta^{-1}$,
   then $\cD=\cD_{14,2}$, $p\not=2, 3$.
   If $q=-1$,
   then $a_{21}^{r_3(X)}=-3$,
   which is a contradiction.

   Case $b_{2b}$.
   Consider the case $q_{11}=q_{22}=-1$, $q_{33}=1$, and $p=3$.
   Then the Dynkin diagram of $r_3(X)$ is
   \Dchainthree{}{$-1$}{$q$}{-$q^{-2}$}{$q$}{$1$}.
   By $a_{23}^{r_3(X)}=-1$ we get $q=-1$ and hence $\cD=\cD_{12',1}$.

   Case $b_{2c}$.
   Consider the case  $q_{11}=q_{22}=1$, $q_{33}:=\zeta\in G_3'$, and $p=2$.
   Then the Dynkin diagram of $r_3(X)$ is
   \Dchainthree{}{$1$}{$q$}{$\zeta q^{-2}$}{$~q{\zeta}^{-1}$}{$\zeta$}~with $\zeta\in G'_3$, $q\not=\zeta$, and $p=2$.
   The condition $a_{23}^{r_3(X)}=-1$ gives $q=\zeta^{-1}$ and hence $\cD=\cD_{52}$, where $q\in G_3'$ and $p=2$.

   Case $b_3$.
   Consider the case $(2)_{q_{11}}=(3)_{q_{33}}=0$, $q_{22}q_1-1=0$, $q_{22}q_2-1=0$, $(2)_{q_{22}}\not=0$, and $(2)_{q_{33}}(q_{33}q_2-1)\not=0$.
   Let $q_{22}:=q$ and $q_{33}:=\zeta$.
   Then $q_1=q_2=q^{-1}$ and $q\not=\zeta$.
   We distinguish subcases $b_{3a}$, $b_{3b}$, and $b_{3c}$.

   Case $b_{3a}$.
   Consider the case $q_{11}=-1$, $\zeta\in G_3'$ and $p\not=2,3$.
   By Lemma~\ref{lem:qij'} we get the Dynkin diagrams of $X$ and $r_3(X)$
  \begin{gather*}
\rule{0pt}{6\unitlength}
\Dchainthree{r}{$-1$}{$q^{-1}$}{$q$}{$q^{-1}$}{$\zeta$} \quad \Rightarrow \quad
\Dchainthree{}{$-1$}{$q^{-1}$}{$\zeta q^{-1}$}{$q{\zeta}^{-1}$}{$\zeta$}\quad \quad
\end{gather*}
  with $\zeta\in G_3'$, $q\notin\{-1, 1, \zeta\}$, and $p\not=2,3$.
  We get $a_{21}^{r_3(X)}\in \{-1, -2\}$ by Definition~\ref{defB3}.
  Hence $((\zeta q^{-1})+1) (\zeta q^{-2}-1)(q^{-3}-1)(\zeta^2q^{-3}-1)=0$.
  Then we obtain $q\in \{{\zeta}^{-1}, -{\zeta}^{-1}, -{\zeta}\}$ or $q^{-3}=\zeta$.
  If $q={\zeta}^{-1}$,
  then  $\cD=\cD_{51}$, where $q\in G'_3$, $p\not=2, 3$.
  If $q=-{\zeta}^{-1}$,
  then  $\cD=\cD_{14,1}$, $p\not=2, 3$.
  If $q=-\zeta$, then by Lemmas~\ref{lem:aij} and~\ref{lem:qij'} we get $A^{r_3(X)}=A^X$ and
  $a_{23}^{r_1r_3(X)}=-3\notin \{-1,-2\}$,
  which is a contradiction.
  If $q^{-3}=\zeta\in G_3'$,
  then $a^{r_1r_3(X)}_{23}=-8\notin \{-1,-2\}$,
  which is again a contradiction.

  Case $b_{3b}$.
  Consider the case $q_{11}=1$, $\zeta\in G_3'$ and $p=2$.
  By Lemma~\ref{lem:qij'} we get the reflection $r_3(X)$ of $X$ is
 \begin{gather*}
\rule{0pt}{5\unitlength}
\Dchainthree{r}{$1$}{$q^{-1}$}{$q$}{$q^{-1}$}{$\zeta$} \quad \Rightarrow \quad
\Dchainthree{}{$1$}{$q^{-1}$}{$\zeta q^{-1}$}{$q{\zeta}^{-1}$}{$\zeta$}\quad \quad
\end{gather*}
  with $\zeta\in G_3'$, $q\not=\zeta$, and $p=2$.
  If $a_{21}^{r_3(X)}=-1$,
  then $q={\zeta}^{-1}$ and $\cD=\cD_{51}$, where $q\in G'_3$, $p=2$.
  If $a_{21}^{r_3(X)}=-2$,
  then $q^{-3}=\zeta\in G_3'$.
  then by Lemmas~\ref{lem:aij} and~\ref{lem:qij'} we get $a^{r_1r_3(X)}_{23}=-8\notin \{-1,-2\}$,
  which is a contradiction.

   Case $b_{3c}$.
  Consider the case $q_{11}=-1$, $\zeta=1$ and $p=3$.
  By Lemma~\ref{lem:qij'} the Dynkin diagrams of $r_3(X)$ is
\Dchainthree{}{$-1$}{$q^{-1}$}{$q^{-1}$}{$q$}{$1$}\quad
  with $q\notin G_2'$, and $p=3$.
  Then $a_{21}^{r_3(X)}\notin \{-1, -2\}$, which is a contradiction.

   Case $b_{4}$.
  Consider the case $(2)_{q_{22}}=(3)_{q_{33}}=0$, $q_{11}q_1-1=0$, $(2)_{q_{11}}\not=0$, and $(2)_{q_{33}}(q_{33}q_2-1)\not=0$.
  Let $q_{11}:=q$ and $q_{33}:=\zeta$.
  Then $q_1=q^{-1}$ and $q_2=q$ by $q_1q_2=1$.
  Hence $q\not=\zeta^{-1}$.
  We distinguish three subcases $b_{4a}$, $b_{4b}$, and $b_{4c}$.

  case $b_{4a}$.
  Consider the case $q_{22}=-1$, $\zeta\in G_3'$, $qq_1-1=0$, and $p\not=2, 3$.
  The Dynkin diagrams of $X$ and $r_3(X)$ are
  \begin{gather*}
\rule{0pt}{5\unitlength}
\Dchainthree{r}{$q$}{$q^{-1}$}{$-1$}{$q$}{$\zeta$} \quad \Rightarrow \quad
\Dchainthree{}{$q$}{$q^{-1}$}{-$\zeta q^2$}{${(q\zeta)}^{-1}$}{$\zeta$}\quad \quad
\end{gather*}
  with $\zeta\in G_3'$, $q\notin \{1, -1, \zeta^{-1}\}$, and $p\not= 2, 3$.
   We get $a_{23}^{r_3(X)}=-1$ by Definition~\ref{defB3}.
   Then $(\zeta q^2-1)(\zeta q^2(\zeta q)^{-1}+1)=0$.
   Hence $q=-\zeta$ or $q=\zeta$.
   If $q=-\zeta$ then $\cD=\cD_{14,3}$, $p\not= 2, 3$.
   If $q=\zeta$
   then $\cD=\cD_{53}$, where $q:=-{\zeta}^{-1}$, $\zeta\in G_3'$, and $p\not= 2, 3$.

   Case $b_{4b}$.
   Consider the case $q_{22}=-1$, $\zeta=1$, and $p=3$.
   Then the Dynkin diagram of $r_3(X)$ is
   \Dchainthree{}{$q$}{$q^{-1}$}{-$q^2$}{${q}^{-1}$}{$1$}~
   with $q\notin G_2'$, $p=3$.
   Then $a_{23}^{r_3(X)}\not=-1$,
   which is a contradiction.

   Case $b_{4c}$.
   Consider the case $q_{22}=1$, $q_{33}:=\zeta\in G_3'$, and $p=2$.
   Then the Dynkin diagram of $r_3(X)$ is
   \Dchainthree{}{$q$}{$q^{-1}$}{$\zeta q^2$}{${(q\zeta)}^{-1}$}{$\zeta$}~.
   The condition $a_{23}^{r_3(X)}=-1$ gives that $q=\zeta$.
   then $\cD=\cD_{53}$, where $q:={\zeta}^{-1}$, $\zeta\in G_3'$, and $p=2$.

   Case $b_5$.
  Consider the case $(2)_{q_{11}}=0$, $q_{22}q_1=1$, $q_{22}q_2=1$, $q_{33}^2q_2=1$, $(2)_{q_{22}}\not=0$, and $(3)_{q_{33}}\not=0$.
  Let $q_{33}:=q$.
  Then $q_2=q^{-2}$, $q_{22}=q^2$ and hence $q_1=q^{-2}$.
  Then $q^2\not=1, -1$, and $q^3\not=1$.
  If $p=3$ and $q_{11}=-1$,
  then $\cD=\cD_{51}$, $p=3$.
  If $p\not=2, 3$ and $q_{11}=-1$,
  then  $\cD=\cD_{51}$, where $q\notin G_3'$ and $p\not=2, 3$.
  If $p=2$ and $q_{11}=1$,
  then $\cD=\cD_{51}$, where $q\notin G_3'$ and $p=2$.

  Case $b_6$.
  Consider the case $(2)_{q_{22}}=0$, $q_{11}q_1=1$, $q_{33}^2q_2=1$, $(2)_{q_{11}}\not=0$, and $(3)_{q_{33}}\not=0$.
  Let $q_{33}:=q$.
  Then $q_2=q^{-2}$ and $q_{1}=q^2$ by $q_1q_2=1$.
  Hence $q_{11}=q^{-2}$.
  If $p=3$ and $q_{22}=-1$,
  then $\cD=\cD_{53}$, $p=3$.
  If $p\not=2, 3$ and $q_{22}=-1$,
  then  $\cD=\cD_{53}$, where $q\notin G_3'$ and $p\not=2, 3$.
  If $p=2$ and $q_{22}=1$,
  then $\cD=\cD_{53}$, where $q\notin G_3'$ and $p=2$.

  Case $b_7$.
  Consider the case $(3)_{q_{33}}=0$, $q_{11}q_1=1$, $q_{22}q_1=1$, $q_{22}q_2=1$, $(2)_{q_{11}}\not=0$, $(2)_{q_{22}}\not=0$ and $(2)_{q_{33}}(q_{33}q_2-1)\not=0$.
  Let $q_{22}:=q$ and $q_{33}:=\zeta$.
  Then $q_2=q_{1}=q^{-1}$, $q_{11}=q$, and $q\not=\zeta$.
  Hence we distinguish two subcases $b_{7a}$ and $b_{7b}$.

  Case $b_{7a}$.
  Consider the case $\zeta\in G_3'$ and $p\not=3$.
  Then
  by Lemma~\ref{lem:qij'} the reflection $r_3(X)$ of $X$ is
  \begin{gather*}
\rule{0pt}{5\unitlength}
\Dchainthree{r}{$q$}{$q^{-1}$}{$q$}{$q^{-1}$}{$\zeta$} \quad \Rightarrow \quad
\Dchainthree{}{$q$}{$q^{-1}$}{$\zeta q^{-1}$}{${q\zeta}^{-1}$}{$\zeta$}\quad \quad
\end{gather*}
  with $\zeta\in G_3'$, $q\notin \{1,-1,\zeta\}$, and $p\not=3$.
  We get $a_{21}^{r_3(X)}\in \{-1, -2\}$ from Definition\ref{defB3}.
  Hence $((\zeta q^{-1})+1) (\zeta q^{-2}-1)(q^{-3}-1)(\zeta^2q^{-3}-1)=0$.
  Then $q\in \{{\zeta}^{-1}, -{\zeta}^{-1}, -{\zeta}\}$ or $q^{-3}=\zeta$.
  If $q={\zeta}^{-1}$,
  then  $\cD=\cD_{21}$, where $q\in G'_3$, $p\not=3$.
  If $q=-{\zeta}^{-1}$ and $p\ne 2$,
  then $\cD=\cD_{12,1}$.
  If $q=-\zeta$ and $p\ne 2$,
  then $\cD=\tau_{321}\cD_{16,5}$.
  If $q^{-3}=\zeta$,
  then $\cD=\cD_{18,1}$.

   Case $b_{7b}$.
  Consider the case $p=3$ and $\zeta=1$.
  Since $(2)_q\not=0$ and $q\not=1$,
  we get $a_{21}^{r_3(X)}\notin \{-1, -2\}$,
  which is a contradiction.

  Case $b_8$.
  Consider the case $q_{11}q_1=q_{22}q_{1}=q_{22}q_{2}=q_{33}^2q_2=1$, $(2)_{q_{11}}\not=0$, $(2)_{q_{22}}\not=0$, and $(3)_{q_{33}}\not=0$.
  If $p\not=2, 3$,
  then  $\cD=\cD_{21}$, where $q\notin G_3'\cup G_4'$ and $p\not=2, 3$.
  If $p=3$,
  then $\cD=\cD_{21}$, where $q\notin G_4'$ and $p=3$.
  If $p=2$,
  then $\cD=\cD_{21}$, where $q\notin G_3'$ and $p=2$.

  Case (c).
  Assume that $A^X$ has a good $C_3$ neighborhood.
  Since $A^X$ is of $C_3$ type,
  we obtain that $q_3=1$, $q_1\not=1$ and $q_2\not=1$ by Lemma~\ref{lem:aij}.
  Since $a_{23}=-2$,
  we get $(2)_{q_{22}}\ne0$ and hence $q_{22}q_1=1$ by $a_{21}=-1$.
  Hence $q_{22}\not=1$.
  Since $a_{23}^{r_1(X)}=-2$,
  we get $(2)_{q_{11}}\not=0$ and hence $q_{11}q_1-1=0$ by $a_{12}=-1$.
  Then we get $q_1q_2=1$ since $a_{13}^{r_2(X)}=0$ by Lemmas~\ref{lem:aij} and~\ref{lem:qij'}.
  Hence,
  by Lemma~\ref{lem:aij} we distinguish the following cases $c_1$, $c_2$, and $c_3$.

  Case $c_{1}$.
  Consider the case $q_{11}q_1=q_{22}q_{1}=q_{33}q_{2}=q_{22}^2q_2=1$.
  Let $q_{22}:=q$.
  Then $q^2\not=1$ and hence $\cD=\cD_{31}$.

  Case $c_2$.
  Consider the case $q_{11}q_1-1=q_{22}q_1-1=q_{22}^2q_2-1=0$, $(2)_{q_{33}}=0$,  and ${q_{33}}q_2-1\not=0$.
  Let $q_{11}:=q$.
  If $p\not=2$,
  then by Lemma~\ref{lem:qij'} we get the reflection of $X$
   \begin{gather*}
\rule{0pt}{5\unitlength}
\Dchainthree{r}{$q$}{$q^{-1}$}{$q$}{$q^{-2}$}{$-1$} \quad \Rightarrow \quad
 \Dchainthree{}{$q$}{$q^{-1}$}{-$q^{-1}$}{$q^2$}{$-1$}\quad \quad
\end{gather*}
with $q\notin G_2'\cup G_4'$.
  We obtain $a_{21}^{r_3(X)}\in \{-1,-2\}$ from Definition~\ref{defC3}.
  Since $q^2\notin G_2'\cup G_4'$,
  we get $a_{21}^{r_3(X)}=-2$ and hence $q^{-3}=-1$ or $q^{-3}=1$.
  If $q^{-3}=-1$, then $p\not=3$ and $\cD=\cD_{13,1}$, with $\zeta \in G_6'$.
  If $q^{-3}=1$, then $p\not=3$ and $\cD=\cD_{13,1}$, with $\zeta \in G_3'$.
  If $p=2$,
  then we get $q^{-3}=1$ by $a_{21}^{r_3(X)}\in \{-1,-2\}$ and then $\cD=\cD_{13',1}$.

  Case $c_3$.
  Consider the case $q_{11}q_1-1=q_{22}q_1-1=0$, $(3)_{q_{22}}=0$, $(2)_{q_{33}}(q_{33}q_2-1)=0$, and $q_{22}^2q_2-1\not=0$.
  Since $(3)_{q_{22}}=0$ and $q_{22}\not=1$,
  we get $q_{22}:=\zeta \in G_3'$, $q_2=\zeta$, and $p\not=3$.
  Hence $q_{22}^2q_2-1=0$, which is a contradiction.
  \end{proof}
%
%
%
%
%
%
%

We point out that by~\cite[Corollary~6]{a-Heck04e} Theorem~\ref{Theo:clasi} yields a classification of three-dimensional braided vector spaces $V$ of diagonal type with finite-dimensional Nichols algebra $\cB(V)$.
\begin{coro}\label{coro-cla}
  Let $\Bbbk$ be a field of characteristic $p>0$.
  Assume that $V$ satisfies the setting in Theorem~\ref{Theo:clasi}.
  Then the Nichols algebra $\cB(V)$ is finite dimensional if and only if
  the Dynkin diagram $\cD$ of $V$ appears in Tables~(\ref{tab.1})-(\ref{tab.3})
  and the labels of the vertices of $\cD$ are roots of unity (including $1$).
  \end{coro}

\setlength{\unitlength}{1mm}
\begin{table}
\centering
\begin{tabular}{r|l|l|l}
row &\text{Dynkin diagrams} & \text{fixed parameters} \\
\hline \hline
 1 & \Dchainthree{}{$q$}{$q^{-1}$}{$q$}{$q^{-1}$}{$q$}
& $q\in \Bbbk^\ast \setminus \{1\}$  \\
\hline
 2 & \Dchainthree{}{$q^2$}{$q^{-2}$}{$q^2$}{$q^{-2}$}{$q$}
& $q\in \Bbbk^\ast \setminus \{1\}$  \\
\hline
 3 & \Dchainthree{}{$q$}{$q^{-1}$}{$q$}{$q^{-2}$}{$q^2$}
& $q\in \Bbbk^\ast \setminus \{1\}$ \\
\hline
 4 & \Dchainthree{}{$1$}{$q^{-1}$}{$q$}{$q^{-1}$}{$q$}\
 \Dchainthree{}{$1$}{$q$}{$1$}{$q^{-1}$}{$q$}
& $q\in \Bbbk^\ast \setminus \{1\}$   \\
\hline
 5 & \Dchainthree{}{$1$}{$q^{-2}$}{$q^2$}{$q^{-2}$}{$q$}\
 \Dchainthree{}{$1$}{$q^2$}{$1$}{$q^{-2}$}{$q$}\
 \Dchainthree{}{$q^2$}{$q^{-2}$}{$1$}{$q^2$}{$q^{-1}$}
& $q\in \Bbbk^\ast \setminus \{1\}$ \\
\hline
 6 & \Dchainthree{}{$1$}{$q^{-1}$}{$q$}{$q^{-2}$}{$q^2$}
& $q\in \Bbbk^\ast \setminus \{1\}$  \\
 & \Dchainthree{}{$1$}{$q$}{$1$}{$q^{-2}$}{$q^2$}\
 \Dtriangle{}{$q$}{$1$}{$1$}{$q^{-1}$}{$q^2$}{$q^{-1}$}
&  \\
\hline
 7 & \Dchainthree{}{$1$}{$q^{-1}$}{$q$}{$q^{-3}$}{$q^3$}
& $q\in \Bbbk^\ast \setminus \{1\}$, $q\notin G'_3$  \\
 & \Dchainthree{}{$1$}{$q$}{$1$}{$q^{-3}$}{$q^3$}\
 \Dtriangle{}{$q$}{$1$}{$1$}{$q^{-1}$}{$q^3$}{$q^{-2}$}\quad
 \Dchainthree{}{$q^3$}{$q^{-3}$}{$1$}{$q^2$}{$q^{-1}$} &\\
\hline
 8 & \Dchainthree{}{$q$}{$q^{-1}$}{$1$}{$q$}{$q^{-1}$}\
 \Dchainthree{}{$1$}{$q^{-1}$}{$q$}{$q^{-1}$}{$1$}
& $q\in \Bbbk^\ast \setminus \{1\}$   \\
 & \Dchainthree{}{$1$}{$q$}{$1$}{$q^{-1}$}{$1$}\
 \Dchainthree{}{$1$}{$q$}{$q^{-1}$}{$q$}{$1$}
&  \\
\hline
 9 & \Dchainthree{}{$q$}{$q^{-1}$}{$1$}{$r^{-1}$}{$r$}
& $q,r,s\in \Bbbk^\ast \setminus \{1\}$, $qrs=1$ \\
& \Dchainthree{}{$q$}{$q^{-1}$}{$-1$}{$s^{-1}$}{$s$}
\Dtriangle{}{$1$}{$1$}{$1$}{$q$}{$r$}{$s$}\quad
 \Dchainthree{}{$r$}{$r^{-1}$}{$1$}{$s^{-1}$}{$s$}&$q\not=r$, $q\not=s$, $r\not=s$ \\
\hline
10 & \Dchainthree{}{$q$}{$q^{-1}$}{$1$}{$q^{-1}$}{$q$}
 & $q\in \Bbbk^\ast \setminus \{1\}$, $q\notin G'_3$ \\
& \Dchainthree{}{$q$}{$q^{-1}$}{$1$}{$q^2$}{$q^{-2}$}
\Dtriangle{}{$1$}{$1$}{$1$}{$q$}{$q$}{$q^{-2}$} &  \\
\hline
11 & \Dchainthree{}{$\zeta $}{$\zeta ^{-1}$}{$1$}{$\zeta ^{-1}$}{$\zeta $}
& $\zeta \in G'_3$ \\
& \rule{24\unitlength}{0pt}\quad \
\Dtriangle{}{$1$}{$1$}{$1$}{$\zeta $}{$\zeta $}{$\zeta $} &\\
\hline
13' & \Dchainthree{}{$\zeta$}{$\zeta^{-1}$}{$\zeta$}{$\zeta $}{$1$}\
   \Dchainthree{}{$\zeta$}{$\zeta^{-1}$}{$\zeta^{-1}$}{$\zeta^{-1}$}{$1$}\
&$\zeta \in G'_3$\\
\hline
15 & \Dchainthree{}{$1$}{$\zeta ^{-1}$}{$\zeta $}{$\zeta $}{$1$} &
$\zeta \in G'_3$ \\
& \Dchainthree{}{$1$}{$\zeta $}{$1$}{$\zeta $}{$1$}\quad \
\Dtriangle{}{$\zeta $}{$1$}{$\zeta $}{$\zeta ^{-1}$}{$\zeta ^{-1}$}{$\zeta ^{-1}$}\quad \
\Dchainthree{}{$1$}{$\zeta ^{-1}$}{$\zeta ^{-1}$}{$\zeta ^{-1}$}{$1$}
&\\
\hline
18 & \Dchainthree{}{$\zeta $}{$\zeta ^{-1}$}{$\zeta $}{$\zeta ^{-1}$}{$\zeta ^{-3}$}\
 \Dchainthree{}{$\zeta $}{$\zeta ^{-1}$}{$\zeta ^{-4}$}{$\zeta ^4$}{$\zeta ^{-3}$}
& $\zeta \in G'_9$ \\
\hline
\end{tabular}
\caption{Dynkin diagrams in characteristic $p=2$}
\label{tab.1}
\end{table}

\setlength{\unitlength}{1mm}
\begin{table}
\centering
\begin{tabular}{r|l|l|l}
row &\text{Dynkin diagrams} & \text{fixed parameters} \\
\hline \hline
 1 & \Dchainthree{}{$q$}{$q^{-1}$}{$q$}{$q^{-1}$}{$q$}
& $q\in \Bbbk^\ast \setminus \{1\}$ \\
\hline
 2 & \Dchainthree{}{$q^2$}{$q^{-2}$}{$q^2$}{$q^{-2}$}{$q$}
& $q\in \Bbbk^\ast \setminus \{-1,1\}$ \\
\hline
 3 & \Dchainthree{}{$q$}{$q^{-1}$}{$q$}{$q^{-2}$}{$q^2$}
& $q\in \Bbbk^\ast \setminus \{-1,1\}$ \\
\hline
 4 & \Dchainthree{}{$-1$}{$q^{-1}$}{$q$}{$q^{-1}$}{$q$}\
 \Dchainthree{}{$-1$}{$q$}{$-1$}{$q^{-1}$}{$q$}
& $q\in \Bbbk^\ast \setminus \{-1,1\}$  \\
\hline
 5 & \Dchainthree{}{$-1$}{$q^{-2}$}{$q^2$}{$q^{-2}$}{$q$}\
 \Dchainthree{}{$-1$}{$q^2$}{$-1$}{$q^{-2}$}{$q$}\
 \Dchainthree{}{$q^2$}{$q^{-2}$}{$-1$}{$q^2$}{$-q^{-1}$}
& $q\in \Bbbk^\ast \setminus \{-1,1\}$, $q\notin G'_4$   \\
\hline
 6 & \Dchainthree{}{$-1$}{$q^{-1}$}{$q$}{$q^{-2}$}{$q^2$}
& $q\in \Bbbk^\ast \setminus \{-1,1\}$   \\
 & \Dchainthree{}{$-1$}{$q$}{$-1$}{$q^{-2}$}{$q^2$}\
 \Dtriangle{}{$q$}{$-1$}{$-1$}{$q^{-1}$}{$q^2$}{$q^{-1}$}
&  \\
\hline
 7 & \Dchainthree{}{$-1$}{$q^{-1}$}{$q$}{$q^{-3}$}{$q^3$}
& $q\in \Bbbk^\ast \setminus \{-1,1\}$  \\
 & \Dchainthree{}{$-1$}{$q$}{$-1$}{$q^{-3}$}{$q^3$}\
 \Dtriangle{}{$q$}{$-1$}{$-1$}{$q^{-1}$}{$q^3$}{$q^{-2}$}\quad \
 \Dchainthree{}{$q^3$}{$q^{-3}$}{$-1$}{$q^2$}{$-q^{-1}$}
 & \\
\hline
 8 & \Dchainthree{}{$q$}{$q^{-1}$}{$-1$}{$q$}{$q^{-1}$}\
 \Dchainthree{}{$-1$}{$q^{-1}$}{$q$}{$q^{-1}$}{$-1$}
& $q\in \Bbbk^\ast \setminus \{-1,1\}$   \\
 & \Dchainthree{}{$-1$}{$q$}{$-1$}{$q^{-1}$}{$-1$}\
 \Dchainthree{}{$-1$}{$q$}{$q^{-1}$}{$q$}{$-1$}
& \\
\hline
 9 & \Dchainthree{}{$q$}{$q^{-1}$}{$-1$}{$r^{-1}$}{$r$}
& $q,r,s\in \Bbbk^\ast \setminus \{1\}$, $qrs=1$ \\
& \Dchainthree{}{$q$}{$q^{-1}$}{$-1$}{$s^{-1}$}{$s$}\
\Dtriangle{}{$-1$}{$-1$}{$-1$}{$q$}{$r$}{$s$}\quad \
\Dchainthree{}{$r$}{$r^{-1}$}{$-1$}{$s^{-1}$}{$s$} &$q\not=r$, $q\not=s$, $r\not=s$  \\
\hline
10 & \Dchainthree{}{$q$}{$q^{-1}$}{$-1$}{$q^{-1}$}{$q$}
 & $q\in \Bbbk^\ast \setminus \{-1,1\}$  \\
& \Dchainthree{}{$q$}{$q^{-1}$}{$-1$}{$q^2$}{$q^{-2}$}
\Dtriangle{}{$-1$}{$-1$}{$-1$}{$q$}{$q$}{$q^{-2}$} &  \\
\hline
12' & \Dchainthree{}{$-1$}{$-1$}{$-1$}{$-1$}{$1$} &\\
\hline
\end{tabular}
\caption{Dynkin diagrams in characteristic $p=3$}
\label{tab.2}
\end{table}
\setlength{\unitlength}{1mm}
\begin{table}
\centering
\begin{tabular}{r|l|l|l}
row &\text{Dynkin diagrams} & \text{fixed parameters}  \\
\hline \hline
 1 & \Dchainthree{}{$q$}{$q^{-1}$}{$q$}{$q^{-1}$}{$q$}
& $q\in \Bbbk^\ast \setminus \{1\}$  \\
\hline
 2 & \Dchainthree{}{$q^2$}{$q^{-2}$}{$q^2$}{$q^{-2}$}{$q$}
& $q\in \Bbbk^\ast \setminus \{-1,1\}$  \\
\hline
 3 & \Dchainthree{}{$q$}{$q^{-1}$}{$q$}{$q^{-2}$}{$q^2$}
& $q\in \Bbbk^\ast \setminus \{-1,1\}$ \\
\hline
 4 & \Dchainthree{}{$-1$}{$q^{-1}$}{$q$}{$q^{-1}$}{$q$}\
 \Dchainthree{}{$-1$}{$q$}{$-1$}{$q^{-1}$}{$q$}
& $q\in \Bbbk^\ast \setminus \{-1,1\}$   \\
\hline
 5 & \Dchainthree{}{$-1$}{$q^{-2}$}{$q^2$}{$q^{-2}$}{$q$}\
 \Dchainthree{}{$-1$}{$q^2$}{$-1$}{$q^{-2}$}{$q$}\
  \Dchainthree{}{$q^2$}{$q^{-2}$}{$-1$}{$q^2$}{$-q^{-1}$}
& $q\in \Bbbk^\ast \setminus \{-1,1\}$  $q\notin G'_4$ \\
\hline
 6 & \Dchainthree{}{$-1$}{$q^{-1}$}{$q$}{$q^{-2}$}{$q^2$}
& $q\in \Bbbk^\ast \setminus \{-1,1\}$ \\
 & \Dchainthree{}{$-1$}{$q$}{$-1$}{$q^{-2}$}{$q^2$}\
 \Dtriangle{}{$q$}{$-1$}{$-1$}{$q^{-1}$}{$q^2$}{$q^{-1}$}
& \\
\hline
 7 & \Dchainthree{}{$-1$}{$q^{-1}$}{$q$}{$q^{-3}$}{$q^3$}
& $q\in \Bbbk^\ast \setminus \{-1,1\}$  \\
 & \Dchainthree{}{$-1$}{$q$}{$-1$}{$q^{-3}$}{$q^3$}\quad \
 \Dtriangle{}{$q$}{$-1$}{$-1$}{$q^{-1}$}{$q^3$}{$q^{-2}$}\
 \Dchainthree{}{$q^3$}{$q^{-3}$}{$-1$}{$q^2$}{$-q^{-1}$}
& $q\notin G'_3$ \\
\hline
 8 & \Dchainthree{}{$q$}{$q^{-1}$}{$-1$}{$q$}{$q^{-1}$}\
\Dchainthree{}{$-1$}{$q^{-1}$}{$q$}{$q^{-1}$}{$-1$}
& $q\in \Bbbk^\ast \setminus \{-1,1\}$   \\
 & \Dchainthree{}{$-1$}{$q$}{$-1$}{$q^{-1}$}{$-1$}\
 \Dchainthree{}{$-1$}{$q$}{$q^{-1}$}{$q$}{$-1$}
& \\
\hline
 9 & \Dchainthree{}{$q$}{$q^{-1}$}{$-1$}{$r^{-1}$}{$r$}
& $q,r,s\in \Bbbk^\ast \setminus \{1\}$, $qrs=1$\\
& \Dchainthree{}{$q$}{$q^{-1}$}{$-1$}{$s^{-1}$}{$s$}\
\Dtriangle{}{$-1$}{$-1$}{$-1$}{$q$}{$r$}{$s$}\quad \
\Dchainthree{}{$r$}{$r^{-1}$}{$-1$}{$s^{-1}$}{$s$} & $q\not=r$, $q\not=s$, $r\not=s$  \\
\hline
10 & \Dchainthree{}{$q$}{$q^{-1}$}{$-1$}{$q^{-1}$}{$q$}
 & $q\in \Bbbk^\ast \setminus \{-1,1\}$  \\
& \Dchainthree{}{$q$}{$q^{-1}$}{$-1$}{$q^2$}{$q^{-2}$}
\Dtriangle{}{$-1$}{$-1$}{$-1$}{$q$}{$q$}{$q^{-2}$} & $q\notin G'_3$  \\
\hline
\end{tabular}
\caption{Dynkin diagrams in characteristic $p>3$}
\end{table}

\addtocounter{table}{-1}
\begin{table}
\centering
\begin{tabular}{r|l|l|l}
row & \text{Dynkin diagrams}  & \text{fixed parameters} \\
\hline \hline
11 & \Dchainthree{}{$\zeta $}{$\zeta ^{-1}$}{$-1$}{$\zeta ^{-1}$}{$\zeta $}
& $\zeta \in G'_3$ \\
& \rule{24\unitlength}{0pt}\
\Dtriangle{}{$-1$}{$-1$}{$-1$}{$\zeta $}{$\zeta $}{$\zeta $} &\\
\hline
12 & \Dchainthree{}{$-\zeta ^{-1}$}{$-\zeta $}{$-\zeta ^{-1}$}{$-\zeta $}{$\zeta $} &
$\zeta \in G'_3$\\
\hline
13 & \Dchainthree{}{$\zeta $}{$\zeta ^{-1}$}{$\zeta $}{$\zeta ^{-2}$}{$-1$}\
 \Dchainthree{}{$\zeta $}{$\zeta ^{-1}$}{$-\zeta^{-1}$}{$\zeta ^2$}{$-1$}
& $\zeta \in G'_3\cup G'_6$ \\
\hline
14 & \Dchainthree{}{$-1$}{$-\zeta $}{$-\zeta ^{-1}$}{$-\zeta $}{$\zeta $}\
 \Dchainthree{}{$-1$}{$-\zeta ^{-1}$}{$-1$}{$-\zeta $}{$\zeta $}\quad \
 \Dchainthree{}{$-\zeta ^{-1}$}{$-\zeta $}{$-1$}{$-\zeta ^{-1}$}{$\zeta ^{-1}$}
& $\zeta \in G'_3$\\
\hline
15 & \Dchainthree{}{$-1$}{$\zeta ^{-1}$}{$\zeta $}{$\zeta $}{$-1$} &
$\zeta \in G'_3$ \\
& \Dchainthree{}{$-1$}{$\zeta $}{$-1$}{$\zeta $}{$-1$}\quad \
\Dtriangle{}{$\zeta $}{$-1$}{$\zeta $}{$\zeta ^{-1}$}{$\zeta ^{-1}$}{$\zeta ^{-1}$}\quad \
\Dchainthree{}{$-1$}{$\zeta ^{-1}$}{$-\zeta ^{-1}$}{$\zeta ^{-1}$}{$-1$}
& \\
\hline
16 & \Dchainthree{}{$-1$}{$\zeta ^{-1}$}{$\zeta $}{$-\zeta ^{-1}$}{$-\zeta $}
 & $\zeta \in G'_3$ \\
 & \Dchainthree{}{$-1$}{$\zeta $}{$-1$}{$-\zeta ^{-1}$}{$-\zeta $}\quad \
\Dtriangle{}{$\zeta $}{$-1$}{$-1$}{$\zeta ^{-1}$}{$-\zeta $}{$-1$} &  \\
 & \Dchainthree{}{$\zeta $}{$-1$}{$-1$}{$-\zeta ^{-1}$}{$-\zeta $}\
\Dchainthree{}{$\zeta $}{$-\zeta ^{-1}$}{$-\zeta $}{$-\zeta ^{-1}$}{$-\zeta $} & \\
\hline
17 & \Dchainthree{}{$-1$}{$-1$}{$-1$}{$\zeta $}{$-1$}
& $\zeta \in G'_3$  \\
 & \Dchainthree{}{$-1$}{$-1$}{$\zeta $}{$\zeta ^{-1}$}{$-1$}\quad \
\Dtriangle{}{$-\zeta $}{$\zeta $}{$-1$}{$-\zeta ^{-1}$}{$\zeta ^{-1}$}{$-\zeta ^{-1}$}\quad \
\Dchainthree{}{$-1$}{$\zeta ^{-1}$}{$\zeta $}{$-\zeta $}{$-1$}&\\
 &\Dchainthree{}{$-1$}{$\zeta $}{$-1$}{$-\zeta $}{$-1$}& \\
 &  \Dchainthree{}{$-1$}{$\zeta $}{$-\zeta $}{$-\zeta ^{-1}$}{$-1$}& \\
 & \Dchainthree{}{$-1$}{$\zeta ^{-1}$}{$\zeta ^{-1}$}{$-\zeta ^{-1}$}{$-1$}\quad \
   \Dtriangle{}{$-1$}{$-1$}{$\zeta $}{$-1$}{$\zeta ^{-1}$}{$-\zeta$}\quad \
   \Dchainthree{}{$-1$}{$-1$}{$-1$}{$-\zeta ^{-1}$}{$\zeta ^{-1}$}&\\
\hline
18 & \Dchainthree{}{$\zeta $}{$\zeta ^{-1}$}{$\zeta $}{$\zeta ^{-1}$}{$\zeta ^{-3}$}\
 \Dchainthree{}{$\zeta $}{$\zeta ^{-1}$}{$\zeta ^{-4}$}{$\zeta ^4$}{$\zeta ^{-3}$}
& $\zeta \in G'_9$  \\
\hline
\end{tabular}
\caption{Dynkin diagrams in characteristic $p>3$}
\label{tab.3}
\end{table}

\setlength{\unitlength}{1mm}
\settowidth{\mpb}{$q_0\in k^\ast \setminus \{-1,1\}$,}
\rule[-3\unitlength]{0pt}{8\unitlength}
\begin{table}
\centering
\begin{tabular}{r|p{7.8cm}|l|l|}
 & \text{exchange graph}  &\text{row in}~\cite[Table~2]{a-Heck04aa} &\text{char} $\Bbbk$\\
\hline \hline
 1 &
 \begin{picture}(2,2)
 \put(1,0){\scriptsize{$\cD_{11}$}}
 \end{picture}
&$1$ & $p>0$\\
\hline
 2 &\begin{picture}(2,4)
 \put(1,0){\scriptsize{$\cD_{21}$}}
 \end{picture}
& $2$ &$p>0$  \\
\hline
 3 &\begin{picture}(2,4)
 \put(1,0){\scriptsize{$\cD_{31}$}}
 \end{picture}
&$3$ &  $p>0$\\
\hline
 4 &
  \begin{picture}(45,4)
 \put(1,0){\scriptsize{$\cD_{41}$}}
 \put(6,0.5){\line(1,0){7}}
 \put(9,1){\scriptsize{$1$}}
 \put(13,0){\scriptsize{$\cD_{42}$}}
 \put(19,0.5){\line(1,0){7}}
 \put(22,1){\scriptsize{$2$}}
 \put(26,0){\scriptsize{$\tau_{321} \cD_{42}$}}
 \put(37,0.5){\line(1,0){6}}
 \put(39,1){\scriptsize{$3$}}
 \put(43,0){\scriptsize{$\tau_{321} \cD_{41}$}}
 \end{picture}
 &$4$ & $p>0$ \\
 \hline
 5 &\begin{picture}(30,3)
 \put(1,0){\scriptsize{$\cD_{51}$}}
 \put(6,0.5){\line(1,0){7}}
 \put(9,1){\scriptsize{$1$}}
 \put(13,0){\scriptsize{$\cD_{52}$}}
  \put(18,0.5){\line(1,0){7}}
 \put(22,1){\scriptsize{$2$}}
 \put(26,0){\scriptsize{$\cD_{53}$}}
 \end{picture}
 & $5$ & $p>0$ \\
 \hline
 6 &\begin{picture}(56,4)
 \put(1,0){\scriptsize{$\cD_{61}$}}
 \put(6,0.5){\line(1,0){7}}
 \put(9,1){\scriptsize{$1$}}
 \put(13,0){\scriptsize{$\cD_{62}$}}
 \put(19,0.5){\line(1,0){7}}
 \put(21,1){\scriptsize{$2$}}
 \put(26,0){\scriptsize{$\cD_{63}$}}
 \put(32,0.5){\line(1,0){6}}
 \put(34,1){\scriptsize{$3$}}
 \put(38,0){\scriptsize{$\tau_{132} \cD_{62}$}}
 \put(49,0.5){\line(1,0){6}}
 \put(51,1){\scriptsize{$1$}}
 \put(55,0){\scriptsize{$\tau_{132} \cD_{61}$}}
 \end{picture}
 &  $6$ &$p>0$\\
 \hline
 7  &\begin{picture}(60,4.5)
 \put(1,0){\scriptsize{$\cD_{71}$}}
 \put(6,1){\line(1,0){7}}
 \put(9,2){\scriptsize{$1$}}
 \put(13,0){\scriptsize{$\cD_{72}$}}
 \put(19,1){\line(1,0){7}}
 \put(21,2){\scriptsize{$2$}}
 \put(26,0){\scriptsize{$\cD_{73}$}}
 \put(32,1){\line(1,0){6}}
 \put(34,2){\scriptsize{$3$}}
 \put(39,0){\scriptsize{$\tau_{231} \cD_{74}$}}
 \end{picture}
 & $7$ &$p>0$ \\
 \hline
 8 &
 \begin{picture}(55,15)
 \put(1,10){\scriptsize{$\cD_{81}$}}
 \put(6,11){\line(1,0){7}}
 \put(9,12){\scriptsize{$2$}}
 \put(13,10){\scriptsize{$\cD_{82}$}}
 \put(18,11){\line(1,0){6}}
 \put(22,12){\scriptsize{$1$}}
 \put(25,10){\scriptsize{$\cD_{83}$}}
 \put(30,11){\line(1,0){6}}
 \put(33,12){\scriptsize{$3$}}
 \put(37,10){\scriptsize{$\tau_{321}\cD_{82}$}}
 \put(35,8){\line(-2,-1){15}}
 \put(25,4){\scriptsize{$1$}}
 \put(49,11){\line(1,0){8}}
 \put(53,12){\scriptsize{$2$}}
 \put(58,10){\scriptsize{$\tau_{321}\cD_{81}$}}
 \put(13,0){\scriptsize{$\cD_{84}$}}
 \put(15,3){\line(0,1){5}}
 \put(13,5){\scriptsize{$3$}}
 \end{picture}
 & $8$ &$p>0$\\
 \hline
9 &\begin{picture}(59,14)
 \put(1,10){\scriptsize{$\cD_{91}$}}
 \put(6,10){\line(1,0){7}}
 \put(9,11){\scriptsize{$2$}}
 \put(13,10){\scriptsize{$\cD_{93}$}}
 \put(18,10){\line(1,0){7}}
 \put(21,11){\scriptsize{$1$}}
 \put(26,10){\scriptsize{$\tau_{213}\cD_{92}$}}
 \put(11,0){\scriptsize{$\tau_{231} \cD_{94}$}}
 \put(15,3){\line(0,1){5}}
 \put(13,4){\scriptsize{$3$}}
 \end{picture}
 & $9$ &$p>0$\\
 \hline
10  &\begin{picture}(59,14)
 \put(1,10){\scriptsize{$\cD_{10,1}$}}
 \put(8,10){\line(1,0){7}}
 \put(11,11){\scriptsize{$2$}}
 \put(15,10){\scriptsize{$\cD_{10,3}$}}
 \put(22,10){\line(1,0){7}}
 \put(25,11){\scriptsize{$1$}}
 \put(29,10){\scriptsize{$\tau_{213}\cD_{10,2}$}}
 \put(11,0){\scriptsize{$\tau_{231} \cD_{10,2}$}}
 \put(16,3){\line(0,1){5}}
 \put(14,4){\scriptsize{$3$}}
 \end{picture}
 & $10$ &$p>0$\\
 \hline
11  &\begin{picture}(59,14)
 \put(1,10){\scriptsize{$\cD_{11,1}$}}
 \put(8,10){\line(1,0){7}}
 \put(11,11){\scriptsize{$2$}}
 \put(15,10){\scriptsize{$\cD_{11,2}$}}
 \put(22,10){\line(1,0){6}}
 \put(25,11){\scriptsize{$1$}}
 \put(29,10){\scriptsize{$\tau_{213}\cD_{11,1}$}}
 \put(13,0){\scriptsize{$\tau_{231} \cD_{11,1}$}}
 \put(16,3){\line(0,1){5}}
 \put(14,5){\scriptsize{$3$}}
 \end{picture}
 & $11$ &$p\not=3$\\
 \hline
12 &\begin{picture}(60,4)
 \put(1,0){\scriptsize{$\cD_{12,1}$}}
\end{picture}
 &$12$ &$p\not=2,3$\\
 \hline
$12'$ &\begin{picture}(60,4)
 \put(1,0){\scriptsize{$\cD_{12',1}$}}
\end{picture}
 &$12$ &$p=3$\\
 \hline
13 &\begin{picture}(60,3)
 \put(1,0){\scriptsize{$\cD_{13,1}$}}
 \put(8,0.5){\line(1,0){8}}
 \put(11,1){\scriptsize{$3$}}
 \put(17,0){\scriptsize{$\cD_{13,2}$}}
\end{picture}
 &$13$ &$p\not=2,3$ \\
 \hline
$13'$ &\begin{picture}(60,3)
 \put(1,0){\scriptsize{$\cD_{13',1}$}}
 \put(9,0.5){\line(1,0){7}}
 \put(12,1){\scriptsize{$2$}}
 \put(18,0){\scriptsize{$\cD_{13',2}$}}
 \end{picture}
 & $13$ &$p=2$ \\
 \hline
14& \begin{picture}(60,4)
 \put(1,0){\scriptsize{$\cD_{14,1}$}}
 \put(8,0.5){\line(1,0){7}}
 \put(11,1){\scriptsize{$1$}}
 \put(15,0){\scriptsize{$\cD_{14,2}$}}
 \put(23,0.5){\line(1,0){7}}
 \put(26,1){\scriptsize{$2$}}
 \put(30,0){\scriptsize{$\cD_{14,3}$}}
\end{picture}
 & $14$ &$p\not=2,3$ \\
 \hline
15  &\begin{picture}(59,14)
 \put(1,10){\scriptsize{$\cD_{15,1}$}}
 \put(8,11){\line(1,0){7}}
 \put(11,11.5){\scriptsize{$1$}}
 \put(15,10){\scriptsize{$\cD_{15,2}$}}
 \put(22,11){\line(1,0){7}}
 \put(25,11.5){\scriptsize{$2$}}
 \put(29,10){\scriptsize{$\cD_{15,3}$}}
 \put(3,4){\line(0,1){5}}
 \put(4,5){\scriptsize{$3$}}
 \put(1,1){\scriptsize{$\cD_{15,4}$}}
 \put(8.5,2){\line(1,0){5}}
 \put(10,2.5){\scriptsize{$1$}}
 \put(14,0){\scriptsize{$\tau_{321} \cD_{15,1}$}}
 \put(17,3){\line(0,1){5}}
 \put(15,5){\scriptsize{$3$}}
 \end{picture}
 & $15$ &$p\not=3$\\
 \hline
16 &
 \begin{picture}(62,21)
 \put(1,18){\scriptsize{$\cD_{16,1}$}}
 \put(9,18){\line(1,0){7}}
 \put(11,18.5){\scriptsize{$1$}}
 \put(17,18){\scriptsize{$\cD_{16,2}$}}
 \put(25,18){\line(1,0){7}}
 \put(29,18.5){\scriptsize{$2$}}
 \put(34,18){\scriptsize{$\cD_{16,4}$}}
 \put(42,18){\line(1,0){7}}
 \put(44,18.5){\scriptsize{$3$}}
 \put(49,18){\scriptsize{$\tau_{132}\cD_{16,3}$}}
 \put(61,18){\line(1,0){7}}
 \put(64,18.5){\scriptsize{$1$}}
 \put(68,18){\scriptsize{$\tau_{132}\cD_{16,5}$}}
 \put(50,12){\line(-2,1){9}}
 \put(47,14){\scriptsize{$1$}}
 \put(46,10){\scriptsize{$\tau_{213}\cD_{16,1}$}}
 \put(37,10){\line(1,0){8}}
 \put(40,11){\scriptsize{$2$}}
 \put(23,10){\scriptsize{$\tau_{213}\cD_{16,2}$}}
 \put(15,10){\line(1,0){7}}
 \put(18,11){\scriptsize{$1$}}
 \put(4,13){\scriptsize{$2$}}
 \put(3,12){\line(0,1){5}}
 \put(1,10){\scriptsize{$\tau_{213}\cD_{16,4}$}}
 \put(24,1){\scriptsize{$\tau_{132}\cD_{16,5}$}}
 \put(3,4){\line(0,1){5}}
 \put(4,6){\scriptsize{$3$}}
 \put(15,1){\line(1,0){8}}
 \put(19,2){\scriptsize{$2$}}
 \put(1,0){\scriptsize{$\tau_{231}\cD_{16,3}$}}
 \end{picture}
 & $16$ &$p\not=2,3$\\
 \hline
17&
 \begin{picture}(59,15)
 \put(1,10){\scriptsize{$\cD_{17,1}$}}
 \put(8,11){\line(1,0){6}}
 \put(11,12){\scriptsize{$3$}}
 \put(15,10){\scriptsize{$\cD_{17,2}$}}
 \put(22,11){\line(1,0){6}}
 \put(25,12){\scriptsize{$2$}}
 \put(30,10){\scriptsize{$\cD_{17,6}$}}
 \put(37,11){\line(1,0){6}}
 \put(40,12){\scriptsize{$3$}}
 \put(44,10){\scriptsize{$\tau_{231}\cD_{17,3}$}}
 \put(57,11){\line(1,0){6}}
 \put(60,12){\scriptsize{$1$}}
 \put(64,10){\scriptsize{$\tau_{231}\cD_{17,4}$}}
 \put(1,0){\scriptsize{$\tau_{213}\cD_{17,9}$}}
 \put(14,0.5){\line(1,0){6}}
 \put(17,1){\scriptsize{$1$}}
 \put(21,0){\scriptsize{$\cD_{17,7}$}}
 \put(28,0.5){\line(1,0){7}}
 \put(31,1){\scriptsize{$3$}}
 \put(21,3.5){\line(-2,1){12}}
 \put(12,4.5){\scriptsize{$2$}}
 \put(36,0){\scriptsize{$\tau_{231}\cD_{17,8}$}}
 \put(42,3){\line(1,1){6}}
 \put(46,5){\scriptsize{$2$}}
 \put(51,1){\line(1,0){8}}
 \put(54,2){\scriptsize{$1$}}
 \put(59,0){\scriptsize{$\tau_{231}\cD_{17,5}$}}
 \put(66,3){\line(0,1){6}}
 \put(67,5){\scriptsize{$2$}}
 \end{picture}
 & $17$ &$p\not=2,3$\\
 \hline
18 &\begin{picture}(17,4)
 \put(1,0){\scriptsize{$\cD_{18,1}$}}
 \put(8,0.5){\line(1,0){7}}
 \put(11,1){\scriptsize{$3$}}
 \put(15,0){\scriptsize{$\cD_{18,2}$}}
 \end{picture}
 & $18$ &$p\not=3$\\
 \hline
\end{tabular}
\caption{The exchange graphs of $\cC(M)$ in Theorem~\ref{Theo:clasi}.}
\label{tab.4}
\end{table}


\newpage
\begin{thebibliography}{10}
\bibitem{inp-Andr14}
Andruskiewitsch,~N.: On finite-dimensional Hopf algebras.
\newblock  Accepted for publication in Proceedings of the International Congress of Mathematicians.
\newblock  arXiv: 1403.7838.

\bibitem{inp-Andr02}
Andruskiewitsch,~N.: About finite dimensional {H}opf algebras.
\newblock \Btxinshort{.}\ {\em Quantum symmetries in theoretical physics and
  mathematics (Bariloche, 2000)\/}, \btxvolumeshort{.}\ 294
  \btxofseriesshort{.}\ {\em Contemp. Math.\/}, \btxpagesshort{.}\ 1--57. Amer.
  Math. Soc. (2002)


\bibitem{a-AndrGr99}
Andruskiewitsch,~N. \btxandshort{.}\ Gra{\~n}a,~M.: Braided {H}opf algebras
  over non-abelian finite groups.
\newblock Bol. Acad. Nac. Cienc. (C{\'o}rdoba).
\newblock \textbf{63}~(1999), 45--78.


\bibitem{a-AHS08}
Andruskiewitsch,~N. \btxandshort{.}\ Heckenberger,~I. and Schneider,~H.-J.:
The {N}ichols algebra of a semisimple {Y}etter-{D}rinfeld module.
\newblock Amer. J. Math.
\newblock \textbf{132}~(2010) no.~6 1493-1547.

\bibitem{a-AndrSchn98}
Andruskiewitsch,~N. \btxandshort{.}\ Schneider,~H.-J.: Lifting of quantum
  linear spaces and pointed {H}opf algebras of order $p^3$.
\newblock J. Algebra.
\newblock \textbf{209}~(1998), 658--691.

\bibitem{a-AndrSchn00}
Andruskiewitsch,~N. \btxandshort{.}\ Schneider,~H.-J.: Finite quantum groups
  and Cartan matrices.
\newblock Adv. Math.
\newblock \textbf{154}~(2000), 1--45.

\bibitem{inp-AndrSchn02}
Andruskiewitsch,~N. \btxandshort{.}\ Schneider,~H.-J.: Pointed {H}opf algebras.
\newblock \Btxinshort{.}\ {\em New {D}irections in {H}opf {A}lgebras\/},
  \btxvolumeshort{.}~43 \btxofseriesshort{.}\ {\em MSRI Publications\/}.
  Cambridge University Press (2002)

\bibitem{a-AndrSchn05}
Andruskiewitsch,~N. \btxandshort{.}\ Schneider,~H.-J.: On the classification of
  finite-dimensional pointed {H}opf algebras.
\newblock Ann. Math.
\newblock \textbf{171}~(2010), 375--417.

\bibitem{Ang1}
Angiono,~I.:
A presentation by generators and relations of Nichols algebras of diagonal
type and convex orders on root systems.
\newblock Accepted for publication in J. Europ. Math. Soc.
\newblock arXiv: 1008.4144.

\bibitem{Ang2}
Angiono, ~I.:
On Nichols algebras of diagonal type.
\newblock J. Reine Angew. Math.
\newblock \textbf{683}~(2013), 189--251.

\bibitem{clw}
Cibils,~C. \btxandshort{.}\  Lauve,~A. \btxandshort{.}\  Witherspoon,~S.:
Hopf quivers and Nichols algebras in positive characteristic.
\newblock Proc. Amer. Math. Soc.
\newblock \textbf{137}~(2009), no.~12, 4029--4041.

%

\bibitem{c-Heck09b}
Cuntz, ~M. \btxandshort{.}\ Heckenberger,~I:
Weyl groupoids with at most three objects.
\newblock Journal of pure and applied algebra
\newblock \textbf{213}~(2009), no.~6, 1112--1128.


\bibitem{c-Heck12a}
Cuntz,~M. \btxandshort{.}\ Heckenberger,~I:
Finite Weyl groupoids of rank three.
\newblock Transactions of the American Mathematical Society
\newblock \textbf{364}~(2012), no.~3, 1369--1393.

\bibitem{c-Heck14a}
Cuntz,~M. \btxandshort{.}\ Heckenberger,~I:
Finite Weyl groupoids.
\newblock Accepted for publication in J. Reine Angew. Math.
\newblock arXiv: 1008.5291.


%
%


\bibitem{a-Heck04d}
Heckenberger,~I.:
Weyl equivalence for rank 2 {N}ichols algebras of diagonal type.
\newblock  Ann.~Uni.~Ferrara,
\newblock \textbf{51(1)}~(2005), 281--289.

\bibitem{a-Heck04aa}
Heckenberger,~I.:
Classification of arithmetic root systems of rank 3.
\newblock Actas del "XVI Coloquio Latinoamericano de Álgebra"
\newblock 227--252(2005)

\bibitem{a-Heck06a}
Heckenberger,~I.:
The {W}eyl groupoid of a {N}ichols algebra of diagonal type,
\newblock Invent. Math.
\newblock \textbf{164}~(2006), 175--188.

\bibitem{a-Heck04bb}
Heckenberger,~I.:
Examples of finite-dimensional rank 2 Nichols algebras of diagonal type.
\newblock Compositio Math.
\newblock \textbf{143(1)}~(2007), 165--190.

\bibitem{a-Heck04e}
Heckenberger,~I.:
Rank 2 Nichols algebras with finite arithmetic root system.
\newblock  Algebra and Representation Theory
\newblock \textbf{11}~(2008), 115--132.


\bibitem{a-Heck09}
Heckenberger,~I.:
Classification of arithmetic root systems.
\newblock Adv.~Math.
\newblock \textbf{220}~(2009), 59-–124.

\bibitem{HS10}
Heckenberger,~I. \btxandshort{.}\ Schneider,~H.-J.:
Root systems and Weyl groupoids for Nichols algebras.
\newblock Proc. Lond. Math. Soc.
\newblock \textbf{101}~(2010), no.~3, 623--654.

\bibitem{HS101}
Heckenberger,~I. \btxandshort{.}\ Schneider,~H.-J.:
Nichols algebras over groups with finite root
system of rank two I.
\newblock J.~Algebra.
\newblock \textbf{324}~(2010), no.~11, 3090--3114.

\bibitem{a-HeckSchn12a}
Heckenberger,~I. \btxandshort{.}\ Schneider,~H.-J.:
Right coideal subalgebras of Nichols algebras and the Duflo order on the Weyl groupoid.
\newblock Israel Journal of Mathematics.
\newblock \textbf{197}~(2013), 139-–187.



\bibitem{VH-14}
Heckenberger,~I.  \btxandshort{.}\ Vendramin,~L.:
A classification of Nichols algebras of semi-simple Yetter-Drinfeld modules over non-abelian groups.
\newblock arXiv: 1412.0857.

\bibitem{WH-14}
Heckenberger,~I.  \btxandshort{.}\ Wang,~J.:
Rank 2 Nichols algebras of diagonal type over fields of positive characteristic.
\newblock SIGMA
\newblock \textbf{11}~(2015), 011, 24 pages

\bibitem{Y-Heck08a}
Heckenberger,~I. \btxandshort{.}\ Yamane,~H.:
A generalization of Coxeter groups, root systems, and Matsumoto's theorem,
\newblock Math.~Z.
\newblock \textbf{259}~(2008), 255--276.




\bibitem{b-Kac90}
Kac,~V.G.:
Infinite dimensional {L}ie algebras.
\newblock Cambridge Univ.~Press
\newblock (1990).

\bibitem{a-Khar99}
Kharchenko,~V.:
A quantum analog of the {P}oincar{\'e}--{B}irkhoff--{W}itt theorem,
\textit{ Algebra and Logic},
\newblock \textbf{38}(4)~(1999), 259--276.


\bibitem{l-2010}
Lusztig,~G.:
Introduction to Quantum Groups.
\newblock Modern Birkh\"auser Classics. Birkh\"auser/Springer, New York, (2010). Reprint of the 1994 edition.


\bibitem{maj05}
Majid,~S.:
Noncommutative differentials and Yang-Mills on permutation groups $S_n$.
\newblock In Hopf algebras in noncommutative geometry and physics, volume 239 of Lecture
Notes in Pure and Appl.~Math.
\newblock ~(2005), 189--213.


\bibitem{n-78}
Nichols,~W.D.:
Bialgebras of type one.
\newblock Commun.~Alg.
\newblock \textbf{6}~(1978), 1521--1552.


\bibitem{Rosso98}
Rosso,~M.:
Quantum groups and quantum shuffles.
\newblock Invent.~Math.
\newblock \textbf{133}~(1998), 399--416.

\bibitem{a-Schauen96}
Schauenburg,~P.:
A characterization of the Borel-like subalgebras of quantum enveloping algebras.
\newblock Commun.~Algebra.
\newblock \textbf{24}~(1996), 2811--2823.

\bibitem{Semi-2011}
Semikhatov,~A.M:
Virasoro central charges for Nichols algebras.
\newblock arXiv: 1109.1767.

\bibitem{Semi-2012}
Semikhatov,~A.M. \btxandshort{.}\ Tipunin,~I.Yu.:
The Nichols algebra of screenings.
\newblock Commun.~Contemp.~Math.
\newblock \textbf{14}~(2012)

\bibitem{Semi-2013}
Semikhatov,~A.M. \btxandshort{.}\ Tipunin,~I.Yu.:
Logarithmic $\hSL2$ CFT{} models from Nichols algebras. I.
\newblock J.~Phys.
\newblock \textbf{A46}~(2013)


\bibitem{w1987}
Woronowicz,~S.L.:
Compact matrix pseudogroups.
\newblock Comm.~Math.~Phys.
\newblock \textbf{111}(4)~(1987), 613–-665.

\bibitem{w1989}
Woronowicz,~S.L.:
Differential calculus on compact matrix pseudogroups (quantum groups).
\newblock Comm.~Math.~Phys.
\newblock \textbf{122}(1)~(1989), 125–-170.
\end{thebibliography}
\end{document}